\UseRawInputEncoding
\documentclass{article}
\pdfoutput=1
\usepackage[a4paper]{geometry}
\usepackage{amssymb,amsmath,amsthm,mathtools,mathrsfs,float}
\usepackage{epstopdf}
\usepackage{multirow}%
\usepackage{amsmath,amssymb,amsfonts}%
\usepackage{amsthm}%
\usepackage{mathrsfs}%
\usepackage[title]{appendix}%
\usepackage{xcolor}%
\usepackage{textcomp}%
\usepackage{manyfoot}%
\usepackage{booktabs}%
\usepackage{algorithm}%
\usepackage{algorithmicx}%
\usepackage{algpseudocode}%
\usepackage[figuresright]{rotating}%
\usepackage{listings}%
\usepackage{cases}%
\usepackage{longtable}
\usepackage{subfig}
\usepackage{lscape}
\usepackage{setspace}
\usepackage{pdftexcmds}
\usepackage{extarrows}

\newtheorem{assumption}{Assumption}[section]
\newtheorem{lemma}{Lemma}[section]
\newtheorem{theorem}{Theorem}[section]
\newtheorem{remark}{Remark}[section]

\title{A Partially Feasible Distributed SQO Method for Two-block General Linearly Constrained Smooth Optimization \footnote{This work was supported by the Natural Science Foundation of China (NFSC) (Grant No. 12171106), the Natural Science Foundation of Guangxi Province (Grant No. 2020GXNSFDA238017), and the NSFC (Grant Nos. 12271113, 12261008).
\newline
\newline
$^1$School of Mathematics and Physics, Guangxi Minzu University,  Center for Applied Mathematics of Guangxi, Nanning 530006, China.
\newline
\newline
$^2$School of Mathematics and Information Science, Guangxi University, Nanning 530004, China.
}}

\author{
Jinbao jian$^1$  \thanks{Email: jianjb@gxu.edu.cn}
\and
Wenrui Chen$^{1, 2}$   \thanks{ Email: wrchen@st.gxu.edu.cn}
\and
Chunming Tang$^{2}$  
\and
Jianghua Yin$^1$ 
}
\date{}

\begin{document}
\maketitle
\begin{abstract}
This paper discusses a class of two-block smooth large-scale optimization problems with both linear equality and linear inequality constraints, which have a wide range of applications, such as
  economic power dispatch, data mining, signal processing, etc.
  Our goal is to develop a novel partially feasible distributed (PFD) sequential
  quadratic optimization (SQO) method (PFD-SQO method) for this kind of problems.
  The design of the method is based on the ideas of SQO method and
  augmented Lagrangian Jacobian splitting scheme
  as well as feasible direction method,
  which decomposes the quadratic optimization (QO) subproblem
  into two small-scale QOs that can be solved independently and parallelly.
  A novel  disturbance contraction term that can be suitably adjusted is
  introduced into the inequality constraints so that the feasible step size
  along the search direction can be increased to 1.
  The new iteration points are generated by the Armijo line search and
  the partially augmented Lagrangian function that only contains equality
  constraints as the merit function. The iteration points always satisfy
  all the inequality constraints of the problem. The theoretical properties
  such as global convergence, iterative complexity,
  superlinear and quadratic rates of convergence of the proposed PFD-SQO
  method are analyzed under appropriate assumptions, respectively.
  Finally, the numerical effectiveness of the method is tested on
  a class of academic examples and an economic power dispatch problem, which shows that the proposed
  method is quite promising.
\end{abstract}
\section{Introduction}\label{section:1}
In this paper, we consider a class of two-block large-scale optimization problems with general linear constraints as follows:
\begin{subequations}\label{0.1}
 \begin{align}
    \min~~ & F(x,y):=f(x)+\theta(y) \label{0.1a}\\
    {\rm s.t.}~~& Ax+By-b=0, \label{0.1b}\\
    & Ex+Fy-d\leqslant 0, \label{0.1c}\\
    & x\in \mathcal{X}:=\{x\in \mathbb{R}^{n_1}|\ l\leqslant Cx\leqslant v\}, \label{0.1d}\\
    & y\in \mathcal{Y}:=\{y\in \mathbb{R}^{n_2}|\ s\leqslant Dy\leqslant r\}, \label{0.1e}
 \end{align}
\end{subequations}
where $f:\mathbb{R}^{n_1}\rightarrow \mathbb{R},
~\theta:\mathbb{R}^{n_2}\rightarrow\mathbb{R}$
are at least first-order continuous differentiable (smooth),
$A\in \mathbb{R}^{m_{1}\times n_{1}}, B\in \mathbb{R}^{m_{1}\times n_{2}}$,
$E\in \mathbb{R}^{m_{2}\times n_{1}}, F\in \mathbb{R}^{m_{2}\times n_{2}}$,
$C\in \mathbb{R}^{l_1\times n_1}, D\in \mathbb{R}^{l_2\times n_2}$,
$l\in \check{\mathbb {R}}^{l_{1}},\ v\in\hat{\mathbb {R}}^{l_{1}}$,
$s\in \check{\mathbb {R}}^{l_{2}}, \ r\in \hat{\mathbb {R}}^{l_{2}}$,
$\check{\mathbb {R}}:=\mathbb{R}\cup\{-\infty\} $ and
$\hat{\mathbb {R}}:=\mathbb{R}\cup\{+\infty\}$.
Call $(n_1, n_2; m_1, m_2; l_1, l_2)$ the scale of problem (\ref{0.1}).
Without loss of generality, we assume that $l<v$ and $s<r$.
Many practical engineering problems can be described as or transformed into the
form of problem (\ref{0.1}), such as compressive sensing \cite{ASS2015},
data mining \cite{WKZXW2003,WLFK2000}, signal processing \cite{YZYW2010},
machine learning \cite{GLY2015}, etc.
For convenience, in this paper, we denote the four special cases of problem (1) as:

 (P1) ---- (\ref{0.1c}) vanishes, and $\mathcal{X}=\mathbb{R}^{n_1}, \mathcal{Y}=\mathbb{R}^{n_2}$;

 (P2) ---- (\ref{0.1c}) vanishes, and $\mathcal{X}$ and $\mathcal{Y}$ are affine manifolds;

 (P3) ---- (\ref{0.1c}) vanishes; 

 (P4) ---- $\mathcal{X}\times \mathcal{Y}=[l,v]\times [p,q]$.

\subsection{ADMM}

The augmented Lagrangian method (ALM) \cite{ALM} is one of
the classical and efficient methods for solving problem (P1), 
and the iterative formula is as follows:
\begin{equation*}\label{ALM}
  \begin{array}{ll}
        (x_{k+1}, y_{k+1})=\arg\min\{ \mathcal{L}_\beta(x,y,\lambda_{k})\mid \ (x, y)\in \mathcal{X}\times \mathcal{Y}\},\\
        \lambda_{k+1}=\lambda_{k}-\beta(Ax_{k+1}+By_{k+1}-b),
        \end{array}
  \end{equation*}
where the augmented Lagrangian function (ALF) 
is defined as
  \begin{equation}\label{La}
    \mathcal{L}_\beta(x,y,\lambda)=f(x)+\theta(y)-
    \lambda^\top (Ax+By-b)+\frac{\beta}{2}\|Ax+By-b\|^2,
  \end{equation}
with $\lambda\in \mathbb{R}^{m_1}$ and $\beta>0$ being the Lagrange multiplier and penalty parameter, respectively.
For small-to-medium-scale problems, ALM usually yields good numerical results.
However, its performance is not satisfactory for large-scale problems.

The distributed optimization methods can decompose a large-scale, challenging optimization problem into several small-scale subproblems that can be solved interactively,
thereby achieving a solution to the original problem.
Alternating Direction Method of Multipliers (ADMM) \cite{GM,GM1,FG,G,He_Yuan_2012,Monteiro_Svaiter_2013,HLWY2014,He_Yuan_2015,HTY2015,BADMM2015,Han_Sun_2015,
Cui_Li_Sun_2015,Fazel_Pong_Sun_Tseng_2013} is one of the effective distributed optimization methods.
ADMM is well suitable for solving large-scale
linearly constrained separable convex optimization problems where
both $f$ and $\theta$ in (\ref{0.1a}) are convex functions.
The classical iterative formula of ADMM for solving problem (P1)
is
\begin{equation*}\label{GS}
 \begin{array}{ll}
      x_{k+1}=\arg\min\limits_{x\in\mathcal{X}}\mathcal{L}_\beta(x,y_{k},\lambda_{k}),\
      y_{k+1}=\arg\min\limits_{y\in\mathcal{Y}}\mathcal{L}_\beta(x_{k+1},y,\lambda_{k}),\\
      \lambda_{k+1}=\lambda_{k}-\beta(Ax_{k+1}+By_{k+1}-b).
 \end{array}
\end{equation*}
The above iteration embeds the Gauss-Seidel scheme \cite{BK}
into each iteration of ALM, which is a serial form of iteration.
Its characteristic is that the subsequent subproblem fully utilizes the information of
the optimal solution generated by the previous subproblems to produce a better solution.
In addition, if the Jacobian scheme \cite{BK} is
embedded into ALM, the iteration formula becomes
\begin{equation*}
\begin{array}{ll}
      x_{k+1}=\arg\min\limits_{x\in\mathcal{X}}\mathcal{L}_\beta(x,y_{k},\lambda_{k}),\
      y_{k+1}=\arg\min\limits_{y\in\mathcal{Y}}\mathcal{L}_\beta(x_{k},y,\lambda_{k}),\\
      \lambda_{k+1}=\lambda_{k}-\beta(Ax_{k+1}+By_{k+1}-b).
\end{array}
\end{equation*}
Obviously, one of the advantages of the Jacobian scheme above is that
the subproblems can be solved in parallel and independently,
reducing the time cost of each iteration.

In recent years, to apply ADMM to more practical problems, both the methods and theories of ADMM and its variants
in the nonconvex case (where at least one of $f$ and $\theta$ in (\ref{0.1a}) is nonconvex) have received much attention \cite{LP2015,HLR2016,BPC,Guo_Han_Wu_2016,WXX_2014,JLMZ2019,TP2020,JLJ,JXC}. 
In the context of the current ADMM for nonconvex optimization,
there are still some issues deserving further investigation.
Firstly, solving two or more low-dimensional subproblems exactly can be as challenging as the original problem
and is still time-consuming, unless the subproblems have very specific structure.
Secondly, for nonlinear constrained or nonconvex optimization problems,
even if the objective function is sufficiently smooth, the corresponding
ADMM is not easy to converge, or the conditions for convergence are
relatively stringent, such as requiring the ALF-based monotone
function to satisfy the K{\L} property.
Thirdly, it is often challenging to achieve superior rates of convergence,
such as superlinear rate of convergence.
Lastly, few ADMMs can effectively handle inequality constraints that involve all variables, even for the linear ones.

\subsection{Distributed SQO method}

The SQO method \cite {Wilson_1963,HR2014,GLR2015,XYZ2015,GKR2016,BCWW2020,BCRZ2021,
Han_1976,Panier,Panier0,Jian_2005,Jian_2006,Jian_2008} has excellent rate of convergence and computational efficiency,
and is one of the most effective methods for solving smooth constrained optimization problems, especially, for small-to-medium-scale problems.
The SQO method, also known as the sequential quadratic programming (SQP) method, is often used to solve the case where the objective functions
$f$ and $\theta$ of problem (P2)
are sufficiently smooth (not necessarily convex).
For the current iteration point $(x_k, y_k)\in \mathcal{X}\times \mathcal{Y}$,
the SQO method usually considers the QO approximation subproblem: 
\begin{equation}\label{QO}
  \begin{array}{ll}
   \min       & f(x_k)\!+\!\nabla f(x_k)^\top(x\!-\!x_k)\!+\!\theta(y_k)\!+\!
   \nabla\theta(y_k)^\top (y\!-\!y_k)\!+\!\frac{1}{2}\parallel (x\!-\!x_k,y\!-\!y_k)\parallel^2_{H_k}\\
   {\rm s.t.} & Ax+By-b=0, x\in \mathcal{X},\ y\in\mathcal{Y},
  \end{array}
\end{equation}
where $\|x\|_H^2:= x^\top H x$, and the matrix $H_k$ is a symmetric positive definite approximation
of the Hessian matrix ${\rm diag}(\nabla^2 f(x_k),\nabla^2\theta(y_k))$ of the complete ALF of (P2).
Obviously, the QO subproblem (\ref{QO}) is the same scale as the original problem.
Therefore, solving (\ref{QO}) is still costly for large-scale problems, even if some efficient solvers are readily available.

Inspired by the splitting idea of ADMM, Jian et al. \cite{Jian_Lao_Chao_Ma_2018}
discussed the case where the constraint sets $\mathcal{X}$ and $\mathcal{Y}$
of problem (P2) 
are non-negative half-spaces.
By choosing the quadratic coefficient matrix as $H_k={\rm diag}(H_k^x,H_k^y)$,
where $H_k^x$ and $H_k^x$ are symmetric positive definite approximations
of $\nabla^2 f(x_k)$ and $\nabla^2\theta(y_k)$, respectively,
the subproblem (\ref{QO}) can be decomposed into two independent QO subproblems by Jacobian splitting:
 \begin{equation*}
 \begin{array}{ll}
 \tilde{x}_{k+1}=\ {\rm arg} \min\limits_{x\geqslant 0} \nabla f(x_k)^\top(x-x_k)
    +\frac{1}{2}\parallel x-x_k\parallel^2_{H^x_k}-\\  \ \ \ \ \ \ \ \ \ \ \ \ \
    \lambda_k^\top(Ax+By_k-b)+\frac{\beta}{2}\|Ax+By_k-b\|^2,
 \end{array}
 \end{equation*}
  \begin{equation*}
 \begin{array}{ll}
 \tilde{y}_{k+1}=\ {\rm arg} \min\limits_{y\geqslant 0} \ \nabla\theta(y_k)^\top (y-y_k)
 +\frac{1}{2}\parallel y-y_k\parallel^2_{H^y_k}-\\  \ \ \ \ \ \ \ \ \ \ \ \ \
 \lambda_k^\top(Ax_k+By-b)
 +\frac{\beta}{2}\|Ax_k+By-b\|^2.
 \end{array}
 \end{equation*}
Under suitable conditions, the global convergence of the method in \cite{Jian_Lao_Chao_Ma_2018}
is proved.

Subsequently, the distributed SQO-type methods were further developed;
see \cite{JOTA,JZY,JZL,JLY,xuxiao}.
In particular, Jian et al. \cite{JZL} considered the two-block nonconvex optimization problem (P3).
For the current feasible iteration point $(x_k,y_k)$,
the following two independent QO subproblems are solved:
\begin{equation*}\label{001.2bx}
\begin{array}{l}
\tilde{x}_{k+1}={\rm arg}\min  \{\nabla f(x^k)^{{\top}}(x-x^k)+\frac{1}{2} \parallel x-x^k\parallel^2_{H^x_k} \ \mid \  A(x-x_k)=0, \ l\leqslant Cx \leqslant v\},\\
\tilde{y}_{k+1}={\rm arg}\min  \{\nabla \theta(y^k)^{{\top}}(y-y^k)+\frac{1}{2} \parallel y-y^k\parallel^2_{H^y_k} \ \mid \  B(y-y_k)=0, \ p\leqslant Dy \leqslant q\}.
\end{array}
\end{equation*}
Then a superlinearly convergent feasible splitting SQO method \cite{JZL} is proposed for the first time.
Due to the fact that the feasible direction generated by the two QO subproblems above has high requirements in terms of equality constraint, that is, the iteration point must satisfies the equality constraint (\ref{0.1b}).
The decline rate of the objective function and the numerical results may be limited to some specific optimizations, i.e., each subproblem is required to have an optimal solution that satisfies the equality constraint. Therefore, the feasible splitting SQO method is worth further investigation in solving the broader optimization problems.

Unlike the problems studied in the related literatures
\cite{Jian_Lao_Chao_Ma_2018,JOTA,JZY,JZL,JLY}, problem
(\ref{0.1}) includes the linear inequality constraint (\ref{0.1c}).
In \cite{xuxiao}, (P4) is discussed.
The main idea is to equivalently transform the linear inequality
constraint (\ref{0.1c}) into $``Ex+Fy+z=d,\ z \geqslant 0 "$ by introducing a relaxation variable $z$,
which further equivalently transforms (P4) into a three-block problem with only linear
equality constraints and box constraint $(x,y)\in [l,v]\times [p,q]$.
A globally convergent ADMM-SQO method \cite{xuxiao} is proposed by combining
the ADMM of the three-block problem with the idea of SQO method.
One of the features of the method is that it has an explicit QO subproblem solution
with respect to the relaxation variable $z$.
However, due to the equivalent transformation, the equality constraint $Ex+Fy+z=d $ is also penalized in ALF, which increases the scale of problem (P4) 
and weakens the feasibility of the obtained approximate solution with respect
to the linear inequality constraint (\ref{0.1c}).

This paper discusses the two-block large-scale smooth optimization problem (\ref{0.1}
with general linear constraints,
aiming to achieve three main goals.
First, the equality constraint (\ref{0.1b}) is still handled by the technique of ALF;
secondly, the inequality constraints (\ref{0.1c})-(\ref{0.1e})
are handled by the feasible direction method;
third, the designed method has fast rate of convergence,
such as superlinear and quadratic rates of convergence.

\subsection{Contributions}
In this paper, a 
partially feasible distributed SQO (PFD-SQO) method is proposed for 
problem (\ref{0.1}). It can be viewed as a hybridization of the SQO method and the distributed ALM as well as the feasible direction method.

$\bullet$ The ALF in this paper only includes the equality
constraint (\ref{0.1b}), which has lower complexity and difficulty of solution compared to
the standard ALF of problem (\ref{0.1}).

$\bullet$  The idea of feasible direction method is used to deal with
the inequality constraints (\ref{0.1c}) - (\ref{0.1e}) so that the iteration points
always satisfy the inequality constraints, i.e., ``partially feasible''.
In particular, in the design of two QO subproblems,
we introduce an adjustable disturbance contraction term
in the inequality constraint (\ref{0.1c}), allowing the feasibility step
size to be increased to a unit step size. This is essential
for the method to have superlinear and quadratic rates of convergence.

$\bullet$  The proposed PFD-SQO method has global convergence under weaker assumption conditions, even when the K{\L} property is not satisfied.
Furthermore, when appropriate second-order approximation conditions are satisfied, the method
achieves superlinear and quadratic rates of convergence, respectively.

$\bullet$ Preliminary numerical experiments on a class of academic examples and an economic power dispatch problem indicate that
the PFD-SQO method is quite promising. Specifically, the PFD-SQO method can efficiently solve the above problems, and the solutions
obtained have good feasibility and optimality, especially with a clear advantage in terms of the time used for calculation.
The numerical results also indicate that the introduction of parameter $c$ is beneficial to the numerical effectiveness of the method.

The main structure of this paper is as follows.
In Section 2, we describe the main idea, iterative steps and basic properties of the method.
In Section 3, the global convergence and iteration complexity are proved under suitable assumptions.
In Section 4, the superlinear and quadratic rates of convergence are demonstrated under appropriate second-order approximation conditions.
In Section 5, we investigate the numerical effectiveness based on two class of mathematical examples and economic power dispatch.

{\bf Notations}: $\mathbb{R}$ denotes the real number set,
$\check{\mathbb {R}}=\mathbb{R}\cup\{-\infty\} $ and
$\hat{\mathbb {R}}=\mathbb{R}\cup\{+\infty\}$;
$(x,y, \cdots):=(x^\top,y^\top,\cdots)^\top$,
where, $x, y, \cdots$ are all column vectors;
component $x_I$ or $(x)_I=(x_i, i\in I)$; $\|x\|^{2}_{A}=x^\top Ax$;
$x\bot y$ denotes $x^\top y=0$;
$\|x\|$ and $\|Q\|$ denote the $\ell_{2}$ norm of the vector and matrix, respectively;
$A\succeq(\succ)0$ means that $A$ is a symmetric semi-positive definite (positive definite) matrix;
$A\succeq(\succ)B$ means that $A-B$ is a symmetric semi-positive definite (positive definite) matrix; $I_n$ means the $n-$order unit matrix;
$[a]$ denotes the largest integer that does not exceed real number $a$.
In order to deal uniformly with the operations on the infinity bound
in the constraints (\ref{0.1d}) and (\ref{0.1e}), we specify the relevant
operations between $\{\pm\infty\}$ and $\mathbb R$ as follows.

\begin{equation*}
  \begin{array}{ll}
\pm\infty+a=\pm\infty, \forall~a\in\mathbb R;
 ~a\times(\pm\infty)=\pm\infty, \forall~a>0;\\
 a\times(\pm\infty)=\mp\infty, \forall~a<0;
 ~a\times(\pm\infty)=0\Leftrightarrow a=0.
  \end{array}
\end{equation*}

\section{Algorithm design}
The (partially) feasible sets of problem (\ref{0.1}) are expressed as
\begin{equation}\label{kexingji}
\begin{array}{ll}
     \mathcal{F}_{\mathfrak{e}}=\{(x, y)|~Ax+By-b=0\},
     ~\mathcal{F}_{\mathfrak{i}}=\{(x, y)|~Ex+Fy-d\leq 0\},\\
     \mathcal{F}_+=(\mathcal{X}\times\mathcal{Y})\cap \mathcal{F}_{\mathfrak{i}},
     ~\mathcal{F}=\mathcal{F}_{\mathfrak{e}}\cap \mathcal{F}_+.
\end{array}
\end{equation}

In this paper, we still consider that the (partially) ALF (\ref{La}) (without considering the constraint
sets $\mathcal{F}_{\mathfrak{i}},\ \mathcal{X}$ and $\mathcal{Y}$) of problem (\ref{0.1}.
In the SQO-type method for problem (\ref{0.1}), the QO subproblem of the current partially feasible iteration point
$(x_k, y_k)\in \mathcal{F}_+$ is considered:
\begin{equation}\label{0.3}
\begin{array}{ll}
\min       & q^{(f,\theta)}_k(x,y):=\nabla f(x_k)^\top(x-x_k)
+\frac{1}{2}\parallel x-x_k\parallel^2_{H^x_k}+  \\
           &\ \ \ \ \ \ \ \ \ \ \ \ \ \ \ \ \ \ \ \ \
\nabla\theta(y_k)^\top (y-y_k)+\frac{1}{2}\parallel y-y_k\parallel^2_{H^y_k}\\
{\rm s.t.} & Ax+By-b=0,\ Ex+Fy-d\leqslant 0,\ x\in \mathcal{X},\ y\in\mathcal{Y}.
\end{array}
\end{equation}
In the subproblem (\ref{0.3}), directly dealing with the equality constraints usually cannot guarantee
the feasibility of the equality constraints, and it will reduce the numerical effectiveness of the method.
So we consider the (partially) ALM of the subproblem (\ref{0.3}):
\begin{equation}\label{0.35}
\begin{array}{ll}
\min       & \mathcal{L}_\beta^
{{\rm SQO}{\text{-}}k}(x,y):=q^{(f,\theta)}_k(x,y)-\lambda_k^\top(Ax+By-b)+\frac{\beta}{2}\parallel Ax+By-b\parallel^2  \\
{\rm s.t.} & Ex+Fy-d\leqslant 0,\ x\in \mathcal{X},\ y\in\mathcal{Y}.
\end{array}
\end{equation}
Obviously, the QO subproblem (\ref{0.35}) is the same scale as problem (\ref{0.1}). Therefore, for large-scale problem,
solving (\ref{0.35}) is still costly.
To overcome the difficulty of scale, absorbing  the ideas of splitting
and reducing dimensions by \cite{Jian_Lao_Chao_Ma_2018,JLY,JOTA,JZL,JZY,xuxiao},
we decompose the subproblem (\ref{0.35}) into two small-scale QO subproblems of Jacobi-type:
\begin{equation*}
\begin{array}{l}
  \min \{\mathcal{L}_\beta^{{\rm SQO}{\text{-}}k}(x,y_k)\ |\ Ex+Fy_k-d\leqslant0,\ x\in \mathcal{X}\},\\
  \min \{\mathcal{L}_\beta^{{\rm SQO}{\text{-}}k}(x_k,y)\ |\ Ex_k+Fy-d\leqslant0,\ y\in \mathcal{Y}\}.
  \end{array}
\end{equation*}
According to (\ref{0.3}) and (\ref{0.35}), two subproblems above can be respectively formulated as
\begin{subequations}\label{0.6}
  \begin{align}
 \min\limits &~~ [\nabla f(x_k)-A^\top(\lambda_k-\beta(Ax_k+By_k-b))]^\top(x-x_k)+
 \frac{1}{2}\|x-x_k\|_{\mathcal{H}_k^x}^2\label{0.6a}\\
 {\rm s.t.} &~~ Ex+Fy_k-d\leqslant 0, x\in \mathcal{X}, \label{0.6b}
  \end{align}
\end{subequations}
and
\begin{subequations}\label{0.7}
  \begin{align}
  \min\limits &~~ [\nabla \theta(y_k)-B^\top(\lambda_k-\beta(Ax_k+By_k-b))]^\top(y-y_k)+
  \frac{1}{2}\|y-y_k\|_{\mathcal{H}_k^y}^2\label{0.7a}\\
 {\rm s.t.} &~~ Ex_k+Fy-d\leqslant 0, y\in \mathcal{Y}, \label{0.7b}
  \end{align}
\end{subequations}
where
\begin{equation}\label{0.701}
\begin{array}{ll}
\mathcal{H}_k^x:=H_k^x+\beta A^\top A\approx\nabla^2f(x_k)+
\beta A^\top A\overset{(\ref{La})}{=} \nabla_{xx}^2\mathcal{L}_\beta(x_k, y_k, \lambda),\\
\mathcal{H}_k^y:=H_k^y+\beta B^\top B\approx\nabla^2\theta(y_k)+\beta B^\top B
\overset{(\ref{La}}{=}\nabla_{yy}^2\mathcal{L}_\beta(x_k, y_k, \lambda).
\end{array}
\end{equation}

To ensure the existence of solutions for the QO subproblems (\ref{0.35}),
(\ref{0.6}) and (\ref{0.7}), the following assumption is required.
\begin{assumption}\label{assumption0}
The following positive definiteness conditions hold:
\begin{subequations}\label{**}
\begin{numcases}{}
\mathcal{H}_k^x=H_k^x+\beta A^\top A\succ 0,
\mathcal{H}_k^y=H_k^y+\beta B^\top B\succ 0, \label{**a}\\
\hat{\mathcal{H}}_k^u:=
\left(
\begin{array}{cc}
\mathcal{H}_k^x & \beta A^\top B \\
\beta B^\top A & \mathcal{H}_k^y
\end{array}
\right)
=\left(
\begin{array}{cc}
H_k^x& \\
 &H_k^y
\end{array}
\right)+\beta(
A~ B)^\top(
A~ B)\succ 0.
\label{**b}
\end{numcases}
\end{subequations}
\end{assumption}

Clearly, Assumption \ref{assumption0} holds whenever $H_k^x\succ 0$ and $H_k^y\succ 0$.
Note that $x_k$ and $y_k$ are feasible solutions of (\ref{0.6}) and (\ref{0.7}), respectively.
If (\ref{**a} holds, the subproblems (\ref{0.6}) and (\ref{0.7}) both have unique optimal
solutions, which are denoted as $\tilde{x}_{k+1}$ and $\tilde{y}_{k+1}$, respectively.
According to (\ref{0.6b}) and (\ref{0.7b}), $\tilde{x}_{k+1}$ and $\tilde{y}_{k+1}$
satisfy $E\tilde{x}_{k+1}+Fy_k-d\leqslant 0$ and $Ex_k+F\tilde{y}_{k+1}-d\leqslant 0$, respectively.
If the search directions are generated by
\begin{equation}\label{0.10}
  d_k^x:=\tilde{x}_{k+1}-x_k, ~d_k^y:=\tilde{y}_{k+1}-y_k,
\end{equation}
one has
$Ed_k^x\leqslant -h_k, Fd_k^y\leqslant -h_k$, where
\begin{equation} \label{hk}
h_k=Ex_k+Fy_k-d\leqslant0.
\end{equation}
Let $t_k$ be the step size and the new iteration point
$(x_{k+1},y_{k+1})=(x_k+t_kd_k^x,~y_k+t_kd_k^y)$. Thus,
\begin{equation*}
Ex_{k+1}+Fy_{k+1}-d=h_k+t_k(Ed_k^x+Fd_k^y)\leqslant(1-2t_k)h_k.
\end{equation*}
It follows from the above relationships that the step size
$t_k\leqslant\frac{1}{2}$ must be controlled to ensure that
$(x_{k+1}, y_{k+1})\in \mathcal{F}_+$.

In order to increase the step size and improve the numerical effect,
we consider further improving the constraints (\ref{0.6b}) and
(\ref{0.7b}), i.e., introducing a freely adjustable parameter $c$
to perturb the right-hand side $0$ of (\ref{0.6b}) and (\ref{0.7b}) to $\frac{c}{2}h_k$.
Then (\ref{0.6}) and (\ref{0.7}) can be improved as
\begin{subequations}\label{0.6X}
  \begin{align}
 \min &~~q_k^f(x):=[\nabla f(x_k)\!-\!A^\top(\lambda_k\!-\!\beta(Ax_k\!+\!By_k\!-\!b))]^\top(x\!-\!x_k)+
 \frac{1}{2}\|x\!-\!x_k\|_{\mathcal{H}_k^x}^2\\
 {\rm s.t.} &~~ Ex+Fy_k-d\leqslant \frac{c}{2}h_k,\\
            &~~ x\in \mathcal{X}=\{x\in \mathbb{R}^{n_1}|\ l\leqslant Cx\leqslant v\},
  \end{align}
\end{subequations}
and
\begin{subequations}\label{0.7Y}
  \begin{align}
 \min &~~q_k^\theta(y):=[\nabla \theta(y_k)\!-vB^\top(\lambda_k\!-\!\beta(Ax_k\!+\!By_k\!-\!b))]^\top(y\!-\!y_k)+
 \frac{1}{2}\|y\!-\!y_k\|_{\mathcal{H}_k^y}^2\\
 {\rm s.t.} &~~ Ex_k+Fy-d\leqslant \frac{c}{2}h_k,\\
            &~~ y\in \mathcal{Y}=\{y\in \mathbb{R}^{n_2}|\ s\leqslant Dy\leqslant r\},
  \end{align}
\end{subequations}
respectively, where the parameter $c\in[0, 1]$.
When $c = 0$, (\ref{0.6X}) and (\ref{0.7Y}) degenerate to
(\ref{0.6}) and (\ref{0.7}) respectively.
Clearly, if Assumption \ref{assumption0}  holds, the QO subproblems
(\ref{0.6X}) and (\ref{0.7Y}) still have unique optimal solutions,
which are still denoted as $\tilde{x}_{k+1}$ and $\tilde{y}_{k+1}$, respectively.
The search direction $(d_k^x,d_k^y)$ is still generated by the
corresponding formula (\ref{0.10}).

For the optimal solutions $\tilde{x}_{k+1}$ and $\tilde{y}_{k+1}$ of
problem (\ref{0.6X}) and (\ref{0.7Y}), by the KKT optimality condition,
there exist corresponding multipliers $\mu_k^x, \mu_k^y \in \mathbb{R}^{m_2}$,\
$\alpha_k^x, \ \gamma_k^x \in \mathbb{R}^{l_1}$ and $ \alpha_k^y,
\gamma_k^y \in \mathbb{R}^{l_2}$, such that
\begin{subequations}\label{0.8X}
  \begin{numcases}{}
\nabla f(x_k)\!+\!\mathcal{H}_k^xd_k^x\!-\!A^\top[\lambda_k\!-\!\beta(Ax_k+By_k-b)]
\!+\!E^\top\mu_k^x\!+\!C^\top(\gamma_k^x-\alpha_k^x)=0, \label{0.8Xa}\\
0\leqslant\mu_k^x\bot (\frac{c}{2}h_k+d-E\tilde{x}_{k+1}-Fy_k)\geqslant0, \label{0.8Xb}\\
0\leqslant\alpha_k^x\bot (C\tilde{x}_{k+1}-l)\geqslant0,\ 0\leqslant\gamma_k^x\bot (v-C\tilde{x}_{k+1})\geqslant0, \label{0.8Xc}
  \end{numcases}
\end{subequations}
and
\begin{subequations}\label{0.9Y}
  \begin{numcases}{}
\nabla \theta(y_k)\!+\!\mathcal{H}_k^yd_k^y\!-\!B^\top[\lambda_k\!-\!\beta(Ax_k+By_k-b)]
\!+\!F^\top\mu_k^y\!+\!D^\top(\gamma_k^y-\alpha_k^y)=0, \label{0.9Ya}\\
0\leqslant\mu_k^y\bot (\frac{c}{2}h_k+d-Ex_k-F\tilde{y}_{k+1})\geqslant0, \label{0.9Yb}\\
0\leqslant\alpha_k^y\bot (D\tilde{y}_{k+1}-s)\geqslant0,\ 0\leqslant\gamma_k^y\bot (r-D\tilde{y}_{k+1})\geqslant0.\label{0.9Yc}
  \end{numcases}
\end{subequations}
\begin{lemma}\label{Proposition1}
  $(\hat{x}, \hat{y})\in \mathbb{R}^{n_1}\times\mathbb{R}^{n_2}$
  is a {\rm KKT} point of problem (\ref{0.1}) if and only if
  there exist corresponding Lagrange  multipliers
  $\hat{\lambda}\in \mathbb{R}^{m_1}$, $\hat{\mu}\in \mathbb{R}^{m_2},
  \hat{\alpha}^x$, $\hat{\gamma}^x\in \mathbb{R}^{l_1}, \hat{\alpha}^y$ and
  $\hat{\gamma}^y\in \mathbb{R}^{l_2}$, such that
  \begin{eqnarray}\label{0.17}
  \left\{\begin{array}{ll}
  &\left(\!
  \begin{array}{c}
  \nabla f(\hat{x}) \\
  \nabla \theta(\hat{y})
  \end{array}
  \!\right)\!-\!
  \left(\!
  \begin{array}{c}
  A^\top \\
  B^\top
  \end{array}
  \!\right)\hat{\lambda}\!+\!
  \left(\!
  \begin{array}{c}
  E^\top \\
  F^\top
  \end{array}
  \!\right)\hat{\mu}\!+\!
  \left(\!
  \begin{array}{cc}
  C^\top&0\\
  0&D^\top
  \end{array}
  \!\right)
  \left(\!
  \begin{array}{c}
  \hat{\gamma}^x-\hat{\alpha}^x \\
  \hat{\gamma}^y-\hat{\alpha}^y
  \end{array}
  \!\right)=
  \left(\!
  \begin{array}{c}
  0 \\
  0
  \end{array}
  \!\right),\\
  &A\hat{x}+B\hat{y}-b=0,\\
  &0\leqslant\hat{\mu}\bot (d-E\hat{x}-F\hat{y})\geqslant0,\\
  &0\leqslant\hat{\alpha}^x\bot (C\hat{x}-l)\geqslant0,\ 0\leqslant\hat{\gamma}^x\bot (v-C\hat{x})\geqslant0,\\
  &0\leqslant\hat{\alpha}^y\bot (D\hat{y}-s)\geqslant0,\ 0\leqslant\hat{\gamma}^y\bot (r-D\hat{y})\geqslant0.
  \end{array}\right.
  \end{eqnarray}
If $(\hat{x}, \hat{y};
\hat{\lambda}, \hat{\mu}, \hat{\alpha}^x,
\hat{\gamma}^x, \hat{\alpha}^y,
\hat{\gamma}^y)$ satisfies (\ref{0.17}), it is said to be a
primal-dual KKT solution or a stationary point of problem (\ref{0.1}).
\end{lemma}

For simplicity, the following notations are used consistently in the rest of the paper:
\begin{equation*}\label{0.155}
  \begin{array}{ll}
w=(x, y, \lambda), ~w_k=(x_k, y_k, \lambda_k), ~u=(x, y),\\
u_k=(x_k, y_k), ~d^u=(d^x, d^y), ~d_k^u=(d_k^x, d_k^y).
  \end{array}
\end{equation*}

From (\ref{0.10}), when $d_k^u=0$ and $Ax_k+By_k-b=0$,
one has $\tilde{x}_{k+1}=x_k$ and $\tilde{y}_{k+1}=y_k$.
Further, from the KKT conditions (\ref{0.8X}), (\ref{0.9Y}) and
(\ref{0.17}), if $\|\mu_k^x-\mu_k^y\|$ is small enough,
the current iteration point $u_k$ can be regarded as an approximate KKT point of problem (\ref{0.1}).
Therefore, inspired by \cite[rule (1.2)]{JZL}, we introduce a recognition criterion
\begin{equation}\label{0.19}
\min\{\|\mu_k^x\|,\ \|\mu_k^y\|\}\leqslant M,~~
\|\mu_k^x-\mu_k^y\|\leqslant M_1(\|d_k^u\|^{\tau_1}+M_2\|Ax_k+By_k-b\|^{\tau_2}),
\end{equation}
where the parameters $\tau_1, \tau_2>0$, $M$ and $M_1$ are appropriately
large, and $M_2\geqslant 0$. The criterion (\ref{0.19}) is used to
test the effectiveness of the QO splitting.
To ensure good convergence of the method, we still consider solving
the QO subproblem (\ref{0.35}) when the QO splitting is invalid, i.e.,
(\ref{0.19}) does not hold.
Since the QO subproblem (\ref{0.35}) always has a feasible solution $u_k$,
it has a unique optimal solution by Assumption \ref{assumption0} (\ref{**b}).
We still denote it by $(\tilde{x}_{k+1}, \tilde{y}_{k+1})$,
and the search direction $d^u_k$ is still generated by (\ref{0.10}),
which corresponds to the KKT condition:
\begin{subequations}\label{0.95}
  \begin{numcases}{}
  \nabla f(x_k)+\mathcal{H}_k^xd^x_k+\beta A^\top Bd^x_k-A^\top[\lambda_k-\beta(Ax_k+By_k-b)]+\nonumber\\
  \ \ \ \  E^\top\mu_k+C^\top(\gamma_k^x-\alpha_k^x)=0, \label{0.95a}\\
  \nabla \theta(y_k)+\mathcal{H}_k^yd^y_k+\beta B^\top Ad^y_k-B^\top[\lambda_k-\beta(Ax_k+By_k-b)]+\nonumber\\
  \ \ \ \  F^\top\mu_k+D^\top(\gamma_k^y-\alpha_k^y)=0, \label{0.95b}\\
  0\leqslant\mu_k\bot (d-E\tilde{x}_{k+1}-F\tilde{y}_{k+1})\geqslant0, \label{0.95c}\\
  0\leqslant\alpha_k^x\bot (C\tilde{x}_{k+1}-l)\geqslant0,~ 0\leqslant\gamma_k^x\bot (v-C\tilde{x}_{k+1})\geqslant0, \label{0.95d}\\
  0\leqslant\alpha_k^y\bot (D\tilde{y}_{k+1}-s)\geqslant0,~ 0\leqslant\gamma_k^y\bot (r-D\tilde{y}_{k+1})\geqslant0. \label{0.95e}
  \end{numcases}
\end{subequations}

Next, we analyze the descent of the ALF $\mathcal{L}_\beta(w)$ of
problem (\ref{0.1}) at the point $w_k$ along the direction $d_k^u$.
It follows from (\ref{La} that
\begin{subequations}\label{nabla_xy}
  \begin{align}
&\nabla_x\mathcal{L}_\beta(w)=\nabla_x\mathcal{L}_\beta(x, y, \lambda)=\nabla f(x)-A^\top[\lambda-\beta(Ax+By-b)], \label{nabla_xya}\\
&\nabla_y\mathcal{L}_\beta(w)=\nabla_y\mathcal{L}_\beta(x, y, \lambda)=\nabla \theta(y)-B^\top[\lambda-\beta(Ax+By-b)]. \label{nabla_xyb}
\end{align}
\end{subequations}
If $(\tilde{x}_{k+1}, \tilde{y}_{k+1})$ is generated by solving the QO subproblems (\ref{0.6X}) and (\ref{0.7Y}),
then by (\ref{0.10}), (\ref{0.8X})-(\ref{0.9Y}),
$u_k\!\in\! \mathcal{X}\!\times\! \mathcal{Y}$ and
$h_k\!\leqslant\!0$, one has
\begin{subequations}\label{xiajiang}
  \begin{align}
 \nabla_x\mathcal{L}_\beta(w_k)^\top d_k^x&=-\|d_k^x\|_{\mathcal{H}_k^x}^2-(Cx_k-l)^\top\alpha_k^x-(v-Cx_k)^\top\gamma_k^x
+(1-\frac{c}{2})h_k^\top\mu_k^x \nonumber\\
&\leqslant -\|d_k^x\|_{\mathcal{H}_k^x}^2,\label{xiajiang1}\\
\nabla_y\mathcal{L}_\beta(w_k)^\top d_k^y&=-\|d_k^y\|_{\mathcal{H}_k^y}^2-(Dy_k-s)^\top\alpha_k^y-(r-Dy_k)^\top\gamma_k^y
+(1-\frac{c}{2})h_k^\top\mu_k^y \nonumber \\
&\leqslant-\|d_k^y\|_{\mathcal{H}_k^y}^2.\label{xiajiang2}
  \end{align}
\end{subequations}
Thus,
\begin{equation}\label{0.145}
\nabla_u\mathcal{L}_\beta(w_k)^\top d_k^u\leqslant-\|d_k^x\|_{\mathcal{H}_k^x}^2-\|d_k^y\|_{\mathcal{H}_k^y}^2=-\|d_k^u\|_{\mathcal{H}_k^u}^2,
\end{equation}
where
\begin{equation}\label{Hku}
\mathcal{H}_k^u=\left(
\begin{array}{cc}
\mathcal{H}_k^x& \\
 &\mathcal{H}_k^y
\end{array}
\right)
=\left(
\begin{array}{cc}
{H}_k^x+\beta A^\top A & \\
 &{H}_k^y+\beta B^\top B
\end{array}
\right).
\end{equation}
When $(\tilde{x}_{k+1}, \tilde{y}_{k+1})$ is generated by solving the QO subproblem (\ref{0.35}), it follows from (\ref{0.95}) that
\begin{equation}\label{0.951}
\begin{array}{ll}
\nabla_u\mathcal{L}_\beta(w_k)^\top d_k^u&=-\|d_k^u\|_{\hat{\mathcal{H}}_k^u}^2\!-\!(Cx_k-l)^\top\alpha_k^x\!-\!
(v-Cx_k)^\top\gamma_k^x-(Dy_k-s)^\top\alpha_k^y-\\
&\ \ \ \ \ (r-Dy_k)^\top\gamma_k^y+h_k^\top\mu_k
\leqslant-\|d_k^u\|_{\hat{\mathcal{H}}_k^u}^2,
\end{array}
\end{equation}
where $\hat{\mathcal{H}}_k^u$ is defined by (\ref{**b}).

If Assumption \ref{assumption0} holds, then by (\ref{0.145}) and (\ref{0.951}),
$\mathcal{L}_\beta(\cdot ,\lambda_k)$ has good descent property along the
direction $d_k^u$ at $u_k$. Based on this, we consider using
$\mathcal{L}_\beta(\cdot ,\lambda_k)$ as the merit function,
starting from $u_k$ along the direction $d_k^u$, generating step size $t_k$ by Armijo line
search, and then generating a new approximate
solution $u_{k+1}=u_k+t_kd_k^u$.

Next, we analyze the maximum feasible step size with respect to the
inequality constraint when moving along the direction $d_k^u$ from
$u_k\in\mathcal{F}_+$.

\begin{lemma}\label{lemma1}
Suppose that parameter $c\in[0,1]$, and $u_k\in \mathcal{F}_+$.
Let $ u_{k+1}(t)=(x_{k+1}(t):=x_k+td_k^x, y_{k+1}(t):=y_k+td_k^y)$.

{\bf(i)} If $(\tilde{x}_{k+1}, \tilde{y}_{k+1})$ is generated by the subproblem
(\ref{0.35}), then for any $t\in[0, 1]$, $u_{k+1}(t)\in \mathcal{F}_+$ holds;

{\bf(ii)} If $\tilde{x}_{k+1}$ and $\tilde{y}_{k+1}$ are generated by the subproblems
(\ref{0.6X}) and (\ref{0.7Y}), respectively, then for any
$t\in [0,\frac{1}{2-c}]$, $u_{k+1}(t)\in \mathcal{F}_+$ holds.
\end{lemma}

{\bf Proof}\ \
(i) From $u_k\in \mathcal{F}_+,\ \tilde{u}_{k+1}:=(
  \tilde{x}_{k+1},\tilde{y}_{k+1})\in \mathcal{F}_+$ and the convexity
  of $\mathcal{F}_+$, one has $u_{k+1}(t)\!=\!(1-t)u_k+t\tilde{u}_{k+1} \in \mathcal{F}_+, \forall~t\in [0,1]$.

(ii) From $u_k, \tilde{u}_{k+1}\in\mathcal{X}\times \mathcal{Y}$
and the convexity of $\mathcal{X}\times \mathcal{Y}$,
one has $u_{k+1}(t)\in\mathcal{X}\times \mathcal{Y}$,
$\forall~ t\in[0,\frac{1}{2-c}]\subseteq [0, 1]$.
Furthermore, $Ed_k^x\leqslant (\frac{c}{2}-1)h_k$ and
$Fd_k^y\leqslant (\frac{c}{2}-1)h_k$ hold from (\ref{0.8X}) and
(\ref{0.9Y}). Therefore
\begin{equation*}
  \begin{array}{ll}
Ex_{k+1}(t)+Fy_{k+1}(t)-d=tEd_k^x+tFd_k^y+h_k\leqslant(2t(\frac{c}{2}-1)+1)h_k\leqslant0,
~\forall~ t\in [0, \frac{1}{2-c}],
  \end{array}
\end{equation*}
i.e., $u_{k+1}(t)\in \mathcal{F}_+$. \hfill$\Box$

\begin{remark}\label{zhu0}
  By Lemma {\rm \ref{lemma1} (ii)},  $d_k^u$ generated by
  the subproblems (\ref{0.6X}) and (\ref{0.7Y}) can reach the longest
  movement $1$ in the partially feasible set $\mathcal{F}_+$
  only if $c=1$. Therefore, the introduction of the contraction term $\frac{c}{2}h_k$
  plays an important role in increasing the line search step size and
  even reaching the unit step size.
\end{remark}

The detailed iterative steps of the {\bf partially feasible distributed
SQO method (abbreviated as PFD-SQO method)} are given below.

\begin{algorithm}
\caption{(PFD-SQO method)}\label{PFD-SQO}
\begin{enumerate}
\setlength{\itemindent}{1em}
\item[\textbf{Step 0}] {\bf (Given initial values)}
Select parameters $\rho\!\in\!(0, 1), \tau_1, \tau_2\!>\!0$,
$c\!\in\![0, 1]$, $\sigma \!\in \!(0, 1)$, $\beta, \xi\!>\!0$,
$M, M_1$ appropriately large, and $M_2\! \geqslant\!0$.
Initial partially feasible point $u_0\!=\!(x_0, y_0)\!\in\! \mathcal{F}_+$,
multiplier $\lambda_0\! \in \! \mathbb{R}^m$, and set $w_0:=\! (u_0, \lambda_0)$.
Select initial $n_1$ and $n_2$ order symmetric positive definite matrices $H_0^x$ and $H_0^y$, set $k\! =\!0$.

\item[\textbf{Step 1}] {\bf (Solve QO subproblems)}
Solve the $x$-QO subproblem (\ref{0.6X}) and the $y$-QO subproblem
(\ref{0.7Y}) in parallel to obtain the optimal solutions
$\tilde{x}_{k+1}$ and $\tilde{y}_{k+1}$, with the corresponding
multipliers $(\mu_k^x, \alpha_k^x, \gamma_k^x)$ and
$(\mu_k^y, \alpha_k^y, \gamma_k^y)$, respectively.

\item[\textbf{Step 2}] {\bf (Generate search direction 1)}
Calculate the search direction $d_k^u=(d_k^x, d_k^y)$ by (\ref{0.10}). If relationship (\ref{0.19}) holds,
set $\tilde{\mathcal{H}}_k^u=\mathcal{H}_k^u$ and
$c_{\max}=\frac{1}{2-c}$, and go to Step 4. Otherwise, go to Step 3.

\item[\textbf{Step 3}] {\bf (Generate search direction 2)}
Solve the $(x, y)$-QO subproblem (\ref{0.35}) to obtain the optimal
solution $(\tilde{x}_{k+1}, \tilde{y}_{k+1})$, with the corresponding
multiplier $(\mu_k, \alpha_k^x, \gamma_x^k, $ $\alpha_k^y, \gamma_k^y)$.
Calculate the search direction $d_k^u$ by (\ref{0.10}),
set $\tilde{\mathcal{H}}_k^u=\hat{\mathcal{H}}_k^u$ and
$c_{\max}=1$, and go to Step 4.

\item[\textbf{Step 4}] {\bf (Armijo line search)}
Calculate the maximum value $t_k\in \{c_{\max}\sigma^i, i=0,1,2,...\}$
that satisfies
\vspace*{-0.2cm}
\begin{equation}\label{0.24}
\mathcal{L}_\beta(u_k+td_k^u,\lambda_k)\leqslant
\mathcal{L}_\beta(w_k)-\rho t\|d_k^u\|_{\tilde{\mathcal{H}}_k^u}^2.
\end{equation}

\item[\textbf{Step 5}] {\bf (Update and optimal identification)}
Generate new iteration point
\begin{subequations}\label{xindedai}
  \begin{numcases}{}
x_{k+1}=x_k+t_k d_k^x,~y_{k+1}=y_k+t_k d_k^y, ~u_{k+1}=(x_{k+1}, y_{k+1}) \label{xindedaia}\\
\lambda_{k+1}=\lambda_k+\xi(Ax_{k+1}+By_{k+1}-b), \label{xindedaib}\\
w_{k+1}=(x_{k+1}, y_{k+1}, \lambda_{k+1}).\label{xindedaic}
\end{numcases}
\end{subequations}
If $(x_{k+1}-x_k, y_{k+1}-y_k, Ax_{k+1}+By_{k+1}-b)=0$,
then $u_{k+1}=u_k$, which is a KKT point of problem (\ref{0.1}), and stop.
Otherwise, compute  new symmetric matrices $H_{k+1}^x$ and $H_{k+1}^y$ such
that they are suitable approximations of $\nabla^2f(x_{k+1})$ and $\nabla ^2\theta(y_{k+1})$, respectively,
and satisfy the positive definiteness requirement (\ref{**}). Let $k:=k+1$, return to Step 1.
\end{enumerate}
\end{algorithm}

From (\ref{La}) and (\ref{xindedaib}), relationship  $\mathcal{L}_\beta(w_{k+1})=
\mathcal{L}_\beta(u_{k+1}, \lambda_k)-\xi\|Ax_{k+1}+By_{k+1}-b\|^2$ holds true.
This, along with the line search (\ref{0.24}), gives that
\begin{equation}\label{0.245}
\mathcal{L}_\beta(w_{k+1})\leqslant \mathcal{L}_\beta(w_k)-\xi\|Ax_{k+1}+By_{k+1}-b\|^2-t_k\rho\|d_k^u\|^2_{\tilde{\mathcal{H}}_k^u}.
\end{equation}
Thus the sequence $\{\mathcal{L}_\beta(w_k)\}$ generated by
Algorithm 1 has excellent monotonic descent property.

To unify the KKT conditions of the QO subproblems in Steps 1
and 3, when $(\tilde{x}_{k+1}, \tilde{y}_{k+1})$  and its KKT multiplier
$(\mu_k, \alpha_k^x, \gamma_k^x, $ $\alpha_k^y,\gamma_k^y)$ are
generated by Step 3, we denote $\mu_k^x=\mu_k^y=\mu_k$.
Thus, from (\ref{0.8X}), (\ref{0.9Y}), (\ref{0.19}), (\ref{0.95}) and
the definition of $\tilde{\mathcal{H}}^u_k$, it follows that whether
$(\tilde{x}_{k+1}, \tilde{y}_{k+1})$ and
$(\mu_k^x, \mu_k^y, \alpha_k^x, \gamma_k^x, \alpha_k^y, \gamma_k^y)$
are generated by Step 1 or Step 3, one always has
\begin{subequations}\label{0.28}
  \begin{numcases}{}
  \left(
  \begin{array}{c}
  \nabla f(x_k) \\
  \nabla \theta(y_k)
  \end{array}
  \right)+\tilde{\mathcal{H}}_k^ud_k^u+
  \left(
  \begin{array}{c}
  0\\
  F^\top(\mu_k^y-\mu_k^x)
  \end{array}
  \right)+\nonumber \\
  \ \ \ \ \left(
  \begin{array}{cccc}
    A^\top & E^\top & C^\top & 0 \\
    B^\top & F^\top & 0 & D^\top
  \end{array}
  \right)
  \left(
  \begin{array}{c}
  \beta(Ax_k+By_k-b)-\lambda_k\\
  \mu_k^x\\
  \gamma_k^x-\alpha_k^x\\
  \gamma_k^y-\alpha_k^y
  \end{array}
  \right)=
  \left(
  \begin{array}{c}
  0 \\
  0
  \end{array}
  \right), \label{0.28a}\\
  0\leqslant\alpha_k^x\bot (C\tilde{x}_{k+1}-l)\geqslant0,~ 0\leqslant\gamma_k^x\bot (v-C\tilde{x}_{k+1})\geqslant0,\label{0.28c}\\
  0\leqslant\alpha_k^y\bot (D\tilde{y}_{k+1}-s)\geqslant0,~ 0\leqslant\gamma_k^y\bot (r-D\tilde{y}_{k+1})\geqslant0,\label{0.28d}\\
  \|\mu_k^x-\mu_k^y\|\leqslant M_1(\!\|d_k^u\|^{\tau_1}+M_2\|Ax_k+By_k-b\|^{\tau_2}\!).\label{0.28f}
  \end{numcases}
\end{subequations}

\begin{lemma}\label{lemma2}
  {\bf(i)} The search direction $d_k^u$ yielded in Step {\rm 4}  satisfies
\begin{equation}\label{0.20}
\nabla_u \mathcal{L}_\beta(w_k)^\top d_k^u\leqslant -\|d_k^u\|_{\tilde{\mathcal{H}}_k^u}^2.
\end{equation}
Therefore, combined with Lemma {\rm \ref{lemma1}}, if $d_k^u\neq0$,
then $d_k^u$ is a feasible descent direction of the constraint
optimization $\min\{\mathcal{L}_\beta (u, \lambda_k)|~u\in \mathcal{F}_+\}$ at point $u_k$
with a maximum feasible step size $c_{\rm max}$.
Further, Step {\rm 4} can be completed after finite computation.
If $d_k^u=0$, then the step size $t_k=c_{\max}\sigma^0=c_{\max}$ is generated by (\ref{0.24}).
Therefore, the line search in Step {\rm 4}  is well-defined, and so is Algorithm 1.

  {\bf(ii)} Algorithm 1 produces a sequence
  $\{u_k\}\subset \mathcal{F}_+$, i.e., $\{u_k\}$ always satisfies
  the inequality constraints of problem (\ref{0.1});

  {\bf(iii)} If $(x_{k+1}-x_k, y_{k+1}-y_k, Ax_{k+1}+By_{k+1}-b)=0$
  holds in Step {\rm 5},
  then $u_{k+1}=u_k$, which is a KKT point of problem (\ref{0.1}).
\end{lemma}

{\bf Proof}\
{(i)} From (\ref{0.145}), (\ref{0.951}) and the composition of the
matrix $\tilde{\mathcal{H}}_k^u$ in Steps 2 and 3, it is known that
the relationship (\ref{0.20}) holds.

(ii) From the generation of  $c_{\rm max}$ in Algorithm 1,
Lemma \ref{lemma1} and the line search in Step 4,
it is clear that $\{u_k\}\subset \mathcal{F}_+$ holds.

{(ii)} Since $(x_{k+1}-x_k, y_{k+1}-y_k, Ax_{k+1}+By_{k+1}-b)=0$,
by (\ref{xindedaia}-(\ref{xindedaib}, $t_k>0$ and (\ref{0.28f}),
one has $d^u_k=0,\ \tilde{x}_{k+1}=x_k, \ \tilde{y}_{k+1}=y_k$,
$0=Ax_{k+1}+By_{k+1}-b=Ax_{k}+By_{k}-b, \ u_{k+1}=u_k\in
\mathcal{F}_{\mathfrak{e}}$ and $\mu_k^x\!=\!\mu_k^y$, denoted as $\mu_k$.
Also by $u_k\! \in \! \mathcal{F}_+$, one has $u_k\in \mathcal{F}$,
i.e., $u_k$ is a feasible solution of problem (\ref{0.1}).
On the other hand, $(\frac{c}{2}-1)\! \in \! [-1,-\frac{1}{2}]$ holds.
It follows from (\ref{0.8Xb}), (\ref{0.9Y}) and (\ref{0.95c}) that $\mu_k\bot(d\! -Ex_k\! -\!Fy_k)$.
Thus, from (\ref{0.28}), $(x_k, y_k, \lambda_k, \mu_k, \alpha_k^x,
\gamma_k^x, \alpha_k^y, \gamma_k^y)$ satisfies (\ref{0.17}).
Further, by Lemma \ref{Proposition1}, $u_{k+1}=u_k$, which is a KKT point
of problem (\ref{0.1}). \hfill$\Box$

\section{Global convergence and iterative complexity 
}

From Step 5 of Algorithm 1 and Lemma \ref{lemma2} (iii),
if Algorithm 1 terminates after a finite number of steps,
a KKT point of problem (\ref{0.1}) is generated.
In this section, we analyze the global convergence of Algorithm 1
in the case of generating an infinite iterative sequence $\{w_k\}$, i.e.,
every accumulating point of the sequence $\{w_k\}$ is
a KKT point of (\ref{0.1}). For this purpose, the following assumption needs
to be guaranteed for  $\{\mathcal{H}_k^u\}$ and $\{\hat{\mathcal{H}}_k^u\}$.
\begin{assumption}\label{assumption1}
In Algorithm 1, the sequences $\{\hat{\mathcal{H}}_k^u\}$
and $\{\mathcal{H}_k^u\}$ of matrices generated by (\ref{**}v and (\ref{Hku}),
respectively, are uniformly positive definite, i.e., there exists a
constant $\eta>0$ such that
\begin{equation}\label{yizhi}
\mathcal{H}_k^u\succeq \eta I_{n_1+n_2}, \ \hat{\mathcal{H}}_k^u\succeq\eta I_{n_1+n_2},
\ \forall~k=0, 1, 2, \ldots.
\end{equation}
Furthermore, if a subsequence $\{w_k\}_K:=\{w_k,k\in K\}$ is bounded,
then the corresponding matrix subsequences $\{\mathcal{H}_k^u\}_K$ and
$\{\hat{\mathcal{H}}_k^u\}_K$ are bounded.
\end{assumption}

Obviously, (\ref{yizhi}) implies
\begin{equation}\label{0.27}
\begin{array}{ll}
\|d^x\|_{\mathcal{H}_k^x}^2 \geqslant\eta \|d^x\|^2,
\|d^y\|_{\mathcal{H}_k^y}^2 \geqslant\eta \|d^y\|^2,
\|d^u\|_{\tilde{\mathcal{H}}_k^u}^2\geqslant\eta\|d^u\|^2,
\ \forall \ k.
\end{array}
\end{equation}
\begin{remark}\label{zhu1}
If $\{H_k^x\}$ and $\{H_k^y\}$ are bounded, then according to (\ref{**})
and (\ref{Hku}), Assumption {\rm \ref{assumption1}} holds when the
matrix $(A~B)$ has column full rank and $\beta$ is sufficiently large,
or when $\{H_k^x\}$ and $\{H_k^y\}$ are uniformly positive.
\end{remark}

Let vectors $c_i, d_i, e_i$ and $f_i$ be the $i{\text{-}}$th rows of matrices
$C, D, E$ and $F$, respectively. And we denote index sets by
$I = \{1, \ldots, l_1\},\ J = \{1, \ldots, l_2\}$,
$S=\{1, \ldots, m_2\}$. For $u=(x,y)\in \mathcal{F}_+$,
the active sets of the inequality constraints of (\ref{0.1}v are denoted
as follows: 
\begin{subequations}\label{0.291}
\begin{numcases}{}
I_x^l=\{i\in I:c_ix = l_i\},~ I^v_x=\{i\in I:c_ix = v_i\},~ I_x =I_x^l\cup I^v_x, \label{0.291a}\\
J^s_y=\{ j\in J:d_j y = s_j\},~ J^r_y =\{ j\in J:d_j y = r_j\},~ J_y= J^s_y\cup J^r_y, \label{0.291b}\\
S_{u}=\{i\in S: e_ix+f_i y-d_i=0\}. \label{0.291c}
\end{numcases}
\end{subequations}

\begin{lemma}\label{lemma3.1}
  Suppose that Assumption {\rm\ref{assumption1}} holds, and the iterative sequence $\{w_k\}$ 
  has a bounded infinite subsequence $\{w_k\}_K$,
  then the set $\Omega$ consisting of all accumulation points of
  $\{w_k\}$ is nonempty. Further, the following conclusions hold.

  {\bf (i)} The corresponding direction sequence $\{d_k^u\}_K$ and the
  subsequence $\{w_{k+1}\}_K:=\{w_{k+1}, k\in K\}$ 
  are both bounded;

  {\bf (ii)}The whole sequence $\{\mathcal{L}_\beta(w_k)\}$ is convergent, and
  \begin{subequations}\label{3.100}
   \begin{align}
  & \lim\limits_{k\rightarrow\infty}\mathcal{L}_\beta(w_k)\!=\!
  \inf\limits_{k}\mathcal{L}_\beta(w_k)
  \!=\!\mathcal{L}_\beta(w_*)\!=\!f(x_*)+\theta(y_*),
  \forall \ w_*=(x_*,y_*,\lambda_*)\in \Omega,
  \label{3.100a}
  \\
  & \lim\limits_{k\rightarrow\infty}(Ax_k+By_k-b)=0,\
  \lim\limits_{k\in K}{d^u_k}=0.\label{3.100b}
  \end{align}
  \end{subequations}

 In addition, 
 each $u_*$ in an accumulation point $w_*:=(u_*,\lambda_*)$ of $\{w_k\}$ is a feasible solution of problem (\ref{0.1}), i.e., $u_*\in \mathcal{F}$.
\end{lemma}

{\bf Proof}\
(i) First, if $(\tilde{x}_{k+1}, \tilde{y}_{k+1})$ is generated by
Step 1, since $x_k$ and $y_k$ are feasible solutions of (\ref{0.6X})
and (\ref{0.7Y}), respectively, then $q^f_k(\tilde{x}_{k+1})\leqslant
q^f_k(x_k)=0$ and $q^\theta_k(\tilde{y}_{k+1})\leqslant q^\theta_k(y_k)=0$.
So, combining (\ref{0.27}), one has
\begin{equation*}
    \begin{array}{ll}
  \eta\|d_k^x\|^2\leqslant\|d_k^x\|_{\mathcal{H}_k^x}^2\leqslant-2\nabla f(x_k)^\top d_k^x+2(d_k^x)^\top A^\top (\lambda_k-\beta(Ax_k+By_k-b)),\\
  \eta\|d_k^y\|^2\leqslant\|d_k^y\|_{\mathcal{H}_k^y}^2\leqslant-2\nabla \theta(y_k)^\top d_k^y+2(d_k^y)^\top B^\top (\lambda_k-\beta(Ax_k+By_k-b)).
  \end{array}
\end{equation*}
These, along with $\|d_k^u\|\leqslant \|d_k^x\|+\|d_k^y\|$, further show that
\begin{equation}\label{lem3.1N1}
  \|d_k^u\|\leqslant 2\eta^{-1}\left(\|\nabla f(x_k)\|+\|\nabla\theta(y_k)\|+(\|A\|+\|B\|)\|\lambda_k-
  \beta(Ax_k+By_k-b)\|\right).
\end{equation}

Second, if $(\tilde{x}_{k+1}, \tilde{y}_{k+1})$ is generated by Step 3,
it follows that $\tilde{q}_k^{(f,\theta)}(\tilde{x}_{k+1},
\tilde{y}_{k+1})\leqslant \tilde{q}_k^{(f,\theta)}(x_k, y_k)=0$.
Combining (\ref{0.27}), one also has
\begin{equation*}
\begin{array}{ll}
\eta\|d_k^u\|^2\leqslant\|d_k^u\|_{\tilde{\mathcal{H}}_k^x}^2 \leqslant &-2\nabla f(x_k)^\top d_k^x
+2(d_k^x)^\top A^\top(\lambda_k-\beta(Ax_k+By_k-b))-\\
&~~2\nabla \theta(y_k)^\top d_k^y+2(d_k^y)^\top B^\top (\lambda_k-\beta(Ax_k+By_k-b)).
\end{array}
\end{equation*}
This further gives
\begin{equation}\label{lem3.1N2}
\|d_k^u\|\leqslant
2\eta^{-1}
\left\|
\begin{array}{c}
\nabla f(x_k)+(\lambda_k-\beta(Ax_k+By_k-b))^\top A\\
\nabla \theta(y_k)+(\lambda_k-\beta(Ax_k+By_k-b))^\top B
\end{array}
\right\|.
\end{equation}
Now, from (\ref{lem3.1N1}), (\ref{lem3.1N2}) and the boundedness of $\{w_k\}_K$,
the boundedness of $\{d_k^u\}_K$ is at hand.
Finally, $\{u_{k+1}\}_K$ is bounded by $ u_{k+1}=u_k+t_k d_k^u$,
and then so is $\{w_{k+1}\}_K$  by (\ref{xindedai}).

(ii) First, it follows from the continuity of $\mathcal{L}_{\beta}(w)$ and
the boundedness of $\{w_k\}_K$ that $\{\mathcal{L}_\beta(w_k)\}_K$ is bounded.
So there exists a positive constant $a$ such that
$|\mathcal{L}_\beta(w_k)|\leqslant a$ holds for any $k\in K$.
Second, for any given non-negative integer $k$,
there exists an associated index $i_k\in K$ such that
$i_k\geqslant k$. This, along with the monotonic descent of
$\{\mathcal{L}_\beta(w_k)\}$, shows that
$\mathcal{L}_\beta(w_k)\geqslant \mathcal{L}_\beta(w_{i_k})\geqslant -a$.
Therefore, the whole sequence $\{\mathcal{L}_\beta(w_k)\}$ is
monotonically decreasing with a lower bound, and hence it is convergent.
Finally, $\lim\limits_{k\rightarrow\infty}\mathcal{L}_\beta(w_k)=
\inf\limits_{k}\mathcal{L}_\beta(w_k)$ holds from the monotonic
descent of $\{\mathcal{L}_\beta(w_k)\}$.

In view of  $\lim\limits_{k\rightarrow\infty}(\mathcal{L}_{\beta}
(w_{k+1})\! -\! \mathcal{L}_{\beta}(w_{k}))\! = \!0$,
taking  the limit of (\ref{0.245}) and combining it with (\ref{0.27}), one has
\begin{equation*}
\begin{array}{ll}
0=\lim\limits_{k\rightarrow\infty}(\mathcal{L}_{\beta}(w_{k+1})-\mathcal{L}_{\beta}(w_{k}))&
\leqslant\ \lim\limits_{k\rightarrow\infty}(-\xi\|Ax_{k+1}+By_{k+1}-b\|^{2}-
t_{k}\rho\|d^{u}_{k}\|^{2}_{\tilde{H}^{u}_{k}})\\
&\leqslant\ \lim\limits_{k\rightarrow\infty}(-\xi\|Ax_{k+1}+By_{k+1}-b\|^{2}-
t_{k}\rho\eta\|d^{u}_{k}\|^{2}).
\end{array}
\end{equation*}
Thus
\begin{equation}\label{3.56}
\lim\limits_{k\rightarrow\infty}t_{k}d^{u}_{k}=0,~
\lim\limits_{k\rightarrow\infty}(Ax_{k+1}+By_{k+1}-b)=0,
\ {\rm i.e.,}\ \lim\limits_{k\rightarrow\infty}(Ax_{k}+By_{k}-b)=0.
\end{equation}

For any given $w_*=(x_*,y_*,\lambda_*)\in \Omega$,
it follows from the second relationship of (\ref{3.56}) that
$Ax_*+By_*-b=0$. This, along with the continuity and the definition
(\ref{La}) of $\mathcal{L}_{\beta}(w)$, shows that $\lim\limits_{k\rightarrow\infty}
\mathcal{L}_\beta(w_k)=\mathcal{L}_\beta(w_*)=f(x_*)+\theta(y_*)$.

In the following, $\lim\limits_{k\in K} d_{k}^{u}=0$ is proved.
If not, there exist $\varepsilon>0$ and $k_{0}\in K$ such that
$\|d_{k}^{u}\|>\varepsilon$ holds for all $k\in K_{0}:=\{k\mid k\in K,~k>k_{0}\}$.
For any $k\in K_{0}$ and sufficiently small positive $t$
(independent of $k$), by Taylor's expansion, formulas (\ref{0.20}),
(\ref{0.27}) and the boundedness of $\{(w_k,d_{k}^{u})\}_{K_{0}}$,
one has
\begin{equation*}
  \begin{array}{ll}
    \mathcal{L}_{\beta}(u_{k}+td_{k},\lambda_{k})&=\mathcal{L}_{\beta}(u_{k},
    \lambda_{k})+t\nabla_{u} \mathcal{L}_{\beta}(u_{k},\lambda_{k})^\top d_{k}^{u}+o(t\|d_{k}^{u}\|)\nonumber\\
    &\leqslant \mathcal{L}_{\beta}(u_{k},\lambda_{k})-t\|d_{k}^{u}\|^{2}_{\tilde{H}_{k}^{u}}+o(t)\nonumber\\
    &=\mathcal{L}_{\beta}(u_{k},\lambda_{k})-t\rho\|d_{k}^{u}\|^{2}_{\tilde{H}_{k}^{u}}
    -t(1-\rho)\|d_{k}^{u}\|^{2}_{\tilde{H}_{k}^{u}}+o(t)\nonumber\\
    &\leqslant \mathcal{L}_{\beta}(u_{k},\lambda_{k})-t\rho\|d_{k}^{u}\|^{2}_{\tilde{H}_{k}^{u}}
    -t\eta(1-\rho)\|d_{k}^{u}\|^{2}+o(t)\nonumber\\
    &\leqslant \mathcal{L}_{\beta}(u_{k},\lambda_{k})-t\rho\|d_{k}^{u}\|^{2}_{\tilde{H}_{k}^{u}}
    -t\eta(1-\rho)\varepsilon^{2}+o(t)\nonumber\\
    &\leqslant \mathcal{L}_{\beta}(u_{k},\lambda_{k})-t\rho\|d_{k}^{u}\|^{2}_{\tilde{H}_{k}^{u}}.\nonumber
\end{array}
\end{equation*}
Thus, it follows from the line search (\ref{0.24}) that
$t_{*}:=\inf\{t_{k}: k\in K_{0}\}>0$. Therefore,
$\lim\limits_{k\in K_{0}}\|t_{k}d_{k}^{u}\|\geqslant t_{*}\varepsilon>0$.
This contradicts $\lim\limits_{k\rightarrow\infty}t_{k}d^{u}_{k}=0$
in (\ref{3.56}), thus $\lim\limits_{k\in K} d_{k}^{u}=0$ is proved.
Futher, result (\ref{3.100b}) together with Lemma {\rm\ref{lemma2} (ii)} shows that $u_*\in \mathcal{F}$.
The proof of Lemma \ref{lemma3.1} is complete. \hfill$\Box$

To ensure that the sequences of KKT multipliers generated by solving
the QO subproblems in Algorithm 1 is reasonably bounded,
and thus further guarantee its global convergence,
the following partial linear independent constraint qualification (PLICQ)
is required for problem (\ref{0.1}).
\begin{assumption}\label{assumption2}
For every accumulation point $w_*=(x_*,y_*,\lambda_*)$ of the
iterative sequence $\{w_k\}$ generated by Algorithm 1,
the inequality constraints (\ref{0.1c})-(\ref{0.1e})
have a row full rank matrix
\begin{equation}\label{A}
  \mathcal{B}_{(x_*, y_*)}:=
  \left(
  \begin{array}{cc}
  E_{S_{u_*}}&F_{S_{u_*}}\\
  C_{I_{x_*}}&0\\
  0&D_{J_{y_*}}\\
  \end{array}
  \right)
  \end{equation}
corresponding to the active constraints at the feasible solution $u_*:=(x_*,y_*)$,
where $E_{S_{u_*}},\ F_{S_{u_*}},\ C_{I_{x_*}}$ and $D_{J_{y_*}}$
are the sub-matrices consisting of the row vectors of the
corresponding index sets of $E,\ F,\ C$ and $D$, respectively.
\end{assumption}

\begin{lemma}\label{lemma3.2}
Suppose that Assumptions {\rm\ref{assumption1}} and {\rm\ref{assumption2}}
hold. If a subsequence $\{w_k\}_{K}$ of $\{w_k\}$ is bounded,
then the corresponding sequence
$\{\Lambda_k:=(\mu^x_k,\mu^y_k, \alpha_k^x, \gamma_k^x,
\alpha_k^y,\gamma_k^y)\}_{K}$ of KKT multipliers is also bounded.
\end{lemma}

{\bf Proof}\  \
By contradiction, suppose that $\{\Lambda_k\}_K$ is unbounded.
Then there exists an infinite subset $\tilde{K}$ of $K$ such that
$\|\Lambda_k\|\stackrel{\tilde{K}}{\longrightarrow}+\infty$.
Notice that $\{w_k\}_{K}$ is bounded and there are only limited ways
to select the active sets in (\ref{0.291}). Without loss of
generality (if necessary, select a subset of $\tilde{K}$), let
\begin{equation*}
\begin{array}{ll}
\lim\limits_{k\in\tilde{K}}w_k=w_*=(x_*, y_*,\lambda_*),\
I_k:= I_{\tilde{x}_{k+1}}\equiv \hat{I},\\
J_k:= J_{\tilde{y}_{k+1}}\equiv \hat{J},
S_k:=S_{(\tilde{x}_{k+1},\tilde{y}_{k+1})}\equiv \hat{S},
~\forall~ k\in \tilde{K}.
\end{array}
\end{equation*}
Thus, by (\ref{0.10}) and Lemma \ref{lemma3.1} (ii),
one has $\lim\limits_{k\in\tilde{K}}(\tilde{x}_{k+1}, \tilde{y}_{k+1})
=\lim\limits_{k\in\tilde{K}} (x _k+d^x_k, y_k+d^y_k)=(x_*,y_*)$.
Therefore, it can be obtained from (\ref{0.291}) that
\begin{equation}\label{Nb10}
I_k\subseteq I_{x_*},\ J_k\subseteq J_{y_*},\ S_k\subseteq S_{u_*},\ \forall\ k\in \tilde{K}.
\end{equation}

A contradiction is derived by two cases as follows.

{\bf Case A.} When $k\in\tilde{K}$ is sufficiently large,
$(\tilde{x}_{k+1}, \tilde{y}_{k+1})$ is generated by Step 1, i.e.,
the splitting iteration is valid.
At this case, by (\ref{0.19}) and
the boundedness of $\{(x_k, y_k, d_k^u)\}_K$, one knows that
$\{(\mu_k^x, \mu_k^y)\}_{\tilde{K}}$ is bounded. On the other hand,
it follows from (\ref{0.8X}) and (\ref{0.9Y}) that
\begin{subequations}\label{0.035}
\begin{numcases}{}
\left(
\begin{array}{c}
\nabla f(x_k) \\
\nabla \theta(y_k)
\end{array}
\right)+\mathcal{H}_k^ud_k^u+
\left(
\begin{array}{c}
E^\top\mu_k^x\\
F^\top\mu_k^y
\end{array}
\right)+\mathcal{C}^\top
\left(
\begin{array}{c}
(\gamma_k^x-\alpha_k^x)_{I_k}\\
(\gamma_k^y-\alpha_k^y)_{J_k}
\end{array}
\right)-\nonumber\\
~~~\left(\begin{array}{c}
A^\top\\
B^\top
\end{array}
\right)(\lambda_k-\beta(Ax_k+By_k-b))=
\left(
\begin{array}{c}
0 \\
0
\end{array}
\right),&\label{0.035a}\\
(\gamma_k^x-\alpha_k^x)_{\bar{I}_k}=0,
(\gamma_k^y-\alpha_k^y)_{\bar{J}_k}=0,\label{0.035b}
\end{numcases}
\end{subequations}
where, $\bar{I}_k\!=\!I\backslash I_k, \bar{J}_k\!=\!J\backslash J_k$,
$\mathcal{C}\!={\rm diag}(C_{I_k}$ and $D_{J_k}).$
From Assumption \ref{assumption2} and (\ref{Nb10}),
one knows that the matrix $\mathcal{C}$ is row full rank.
Then by (\ref{0.035a}), one has
\begin{equation*}\label{0.36}
\begin{array}{ll}
\left(
\begin{array}{c}
(\gamma_k^x-\alpha_k^x)_{I_k}\\
(\gamma_k^y-\alpha_k^y)_{J_k}
\end{array}
\right)=&-(\mathcal{C}\mathcal{C}^\top)^{-1}\mathcal{C}\bigg[
\left(
\begin{array}{c}
\nabla f(x_k) \\
\nabla \theta(y_k)
\end{array}
\right)+\mathcal{H}_k^ud_k^u+
\left(
\begin{array}{c}
E^\top\mu_k^x\\
F^\top\mu_k^y
\end{array}
\right)-\\
&~~\left(
\begin{array}{c}
A^\top\\
B^\top
\end{array}
\right)(\lambda_k-\beta(Ax_k+By_k-b))\bigg].
\end{array}
\end{equation*}
This, along with the boundedness of $\{(w_k, d_k^u, \mu_k^x, \mu_k^y,
\mathcal{H}_k^u)\}_K$ and (\ref{0.035b}),
shows that the sequence $\{(\gamma_k^x-\alpha_k^x,
\gamma_k^y-\alpha_k^y)\}_{\tilde{K}}$ is bounded.
Again, it follows from (\ref{0.28c}) and (\ref{0.28d}) that
$(\gamma_k^x)^\top\alpha_k^x=0$ and $(\gamma_k^y)^\top\alpha_k^y=0$, respectively.
Thus,
\begin{equation*}
  \begin{array}{ll}
&\|\alpha_k^x\|^2=(\alpha_k^x)^\top(\alpha_k^x-\gamma_k^x)\leqslant
\|\alpha_k^x\|\cdot\|\alpha_k^x-\gamma_k^x\|\Rightarrow \|\alpha_k^x\|\leqslant\|\gamma_k^x-\alpha_k^x\|,\\
&\|\alpha_k^y\|^2=(\alpha_k^y)^\top(\alpha_k^y-\gamma_k^y)\leqslant
\|\alpha_k^y\|\cdot\|\alpha_k^y-\gamma_k^y\|\Rightarrow \|\alpha_k^y\|\leqslant\|\gamma_k^y-\alpha_k^y\|.
  \end{array}
\end{equation*}
Therefore, $\{(\alpha_k^x,\gamma_k^x, \alpha_k^y, \gamma_k^y)\}_{\tilde{K}}$
is bounded. Thus, the subsequence $\{\Lambda_k\}_{\tilde{K}}$
is bounded, which contradicts the fact that $\|\Lambda_k\|\stackrel{\tilde{K}}{\longrightarrow}+\infty$.

{\bf Case B.} There exists an infinite subsequence
$\tilde{K}^-\subseteq \tilde{K}$ such that for each $k\in \tilde{K}^-$,
$(\tilde{x}_{k+1}, \tilde{y}_{k+1})$ is generated by Step 3, i.e.,
the splitting iteration fails. At this case, one has from (\ref{0.95})
\begin{subequations}\label{0.035B}
\begin{numcases}{}
\left(
\begin{array}{c}
\nabla f(x_k) \\
\nabla \theta(y_k)
\end{array}
\right)+\tilde{\mathcal{H}}_k^ud_k^u+
\mathcal{B}^\top
\left(
\begin{array}{c}
(\mu_k)_{S_k}\\
(\gamma_k^x-\alpha_k^x)_{I_k}\\
(\gamma_k^y-\alpha_k^y)_{J_k}
\end{array}
\right)- \nonumber\\
~~~~~\left(
\begin{array}{c}
A^\top\\
B^\top
\end{array}
\right)(\lambda_k-\beta(Ax_k+By_k-b))=
\left(
\begin{array}{c}
0 \\
0
\end{array}
\right),\label{0.035Ba}\\
(\gamma_k^x-\alpha_k^x)_{\bar{I}_k}=0,
(\gamma_k^y-\alpha_k^y)_{\bar{J}_k}=0,
(\mu_k)_{\bar{S}_k}=0,\label{0.035Bb}
\end{numcases}
\end{subequations}
where $\bar{S}_k=S\setminus S_k$ and
$\mathcal{B}=\mathcal{B}_{(\tilde{x}_{k+1},\tilde{y}_{k+1})}$.
By relationship (\ref{Nb10}) and Assumption \ref{assumption2},
the matrix $\mathcal{B}$ has full rank. Thus, one has from (\ref{0.035Ba})
\begin{equation*}\label{0.36B}
\begin{array}{ll}
\left(
\begin{array}{c}
(\mu_k)_{S_k}\\
(\gamma_k^x-\alpha_k^x)_{I_k}\\
(\gamma_k^y-\alpha_k^y)_{J_k}
\end{array}
\right)=&-(\mathcal{B}\mathcal{B}^\top)^{-1}\mathcal{B}\bigg[
\left(
\begin{array}{c}
\nabla f(x_k) \\
\nabla \theta(y_k)
\end{array}
\right)+\tilde{\mathcal{H}}_k^ud_k^u- \\
&~~\left(
\begin{array}{c}
A^\top\\
B^\top
\end{array}
\right)(\lambda_k-\beta(Ax_k+By_k-b))\bigg].
\end{array}
\end{equation*}
This, together with the boundedness of $\{(w_k, d_k^u,
\tilde{\mathcal{H}}_k^u)\}_K$ and (\ref{0.035Bb}),
shows that the sequence $\{(\mu_k, \gamma_k^x-\alpha_k^x,
\gamma_k^y-\alpha_k^y)\}_{\tilde{K}^-}$ is bounded.
Furthermore, similarly to case A, it can be proven that
$\{(\alpha_k^x,\gamma_k^x, \alpha_k^y, \gamma_k^y)\}_{\tilde{K}^-}$
is bounded. Therefore, the subsequence $\{\Lambda_k\}_{\tilde{K}^-}$
is bounded, which also contradicts $\|\Lambda_k\|\stackrel{\tilde{K}}\longrightarrow+\infty$.
The proof is completed. \hfill$\Box$

Based on the preparation above, the global convergence of Algorithm 1 is established below.

\begin{theorem}\label{theorem1}

Suppose that Assumptions {\rm \ref{assumption1}} and
{\rm \ref{assumption2}} hold. Then, for every accumulation point
$w_*:=(x_*,y_*,\lambda_*)$ of sequence $\{w_k\}$ generated by Algorithm 1,
there exists a subsequence
$\{(\mu_k^x, \mu_k^y, \alpha_k^x, \gamma_k^x, \alpha_k^y, \gamma_k^y)\}_K$ of multipliers
such that $(\mu_k^x, \mu_k^y, \alpha_k^x, \gamma_k^x, \alpha_k^y, \gamma_k^y)\stackrel{K}{\longrightarrow}
(\mu_*, \mu_*,  \alpha_*^x, \gamma_*^x, \gamma_*^y, \gamma_*^y)$ and
$(w_*, \mu_*,  \alpha_*^x, \gamma_*^x, \gamma_*^y, \gamma_*^y)$
constitutes a primal-dual {\rm KKT} solution of problem (\ref{0.1}).
Therefore, Algorithm 1 is globally convergent in this sense.
\end{theorem}

{\bf Proof}\ \
In view of Lemmas
\ref{lemma3.1} and \ref{lemma3.2}, there exists an infinite set
$K$ of iterative indices such that
\begin{subequations}\label{0.29}
\begin{numcases}{}
w_k\stackrel{K}{\longrightarrow} w_*,\ d_k^u 
\stackrel{K}{\longrightarrow} 0, ~u_{k+1}
=u^k+t_kd^u_k\stackrel{K}{\longrightarrow} u_*:=(x_*,y_*), \label{0.29aa} \\
\tilde{u}_{k+1}=:(\tilde{x}_{k+1},\tilde{y}_{k+1})=u^k+d^u_k\stackrel{K}{\longrightarrow} u_*,\label{0.29a}\\
(\mu_k^x, \mu_k^y, \alpha_k^x, \gamma_k^x, \alpha_k^y, \gamma_k^y)\stackrel{K}{\longrightarrow} ( \mu_*^x, \mu_*^y,  \alpha_*^x, \gamma_*^x, \gamma_*^y, \gamma_*^y).
\end{numcases}
\end{subequations}
Thus, it follows from (\ref{0.28f}), (\ref{3.100b}) and (\ref{0.29}) that $\mu_*^x=\mu_*^y:=\mu_*$.
Further, we take $k\stackrel{K}{\longrightarrow}\infty$ in (\ref{0.28})
and combine it with (\ref {0.29}) to obtain
\begin{equation}\label{0.347}
\left\{\begin{array}{ll}
\left(\!
\begin{array}{c}
\nabla f(x_*) \\
\nabla \theta(y_*)
\end{array}
\!\right)\!-\!
\left(\!
\begin{array}{c}
A^\top\\
B^\top
\end{array}
\!\right)\lambda_*\!+\!
\left(\!
\begin{array}{c}
E^\top\\
F^\top
\end{array}
\!\right)\mu_*\!+\!
\left(\!
\begin{array}{cc}
C^\top & 0\\
0 & D^\top
\end{array}
\!\right)\!\left(\!
\begin{array}{c}
\gamma_*^x-\alpha_*^x\\
\gamma_*^y-\alpha_*^y
\end{array}
\!\right)\!=\!
\left(\!
\begin{array}{c}
0 \\
0
\end{array}
\!\right), \\
0\leqslant\alpha_*^x\bot (C{x}_*-l)\geqslant0, ~~0\leqslant\gamma_*^x\bot (v-C{x}_*)\geqslant0,\\
0\leqslant\alpha_*^y\bot (D{y}_*-s)\geqslant0, ~~0\leqslant\gamma_*^y\bot (r-D{y}_*)\geqslant0.
\end{array}\right.
\end{equation}
In addition, notice that
\begin{equation*}
\begin{array}{ll}
\frac{c}{2}h_k+d-E\tilde{x}_{k+1}-Fy_k\stackrel{K}{\longrightarrow}(1-\frac{c}{2})(d-Ex_*-Fy_*),\\
\frac{c}{2}h_k+d-Ex_k-F\tilde{y}_{k+1}\stackrel{K}{\longrightarrow}(1-\frac{c}{2})(d-Ex_*-Fy_*),\\
d-E\tilde{x}_{k+1}-F\tilde{y}_{k+1}\stackrel{K}{\longrightarrow}d-Ex_*-Fy_*.\\
\end{array}
\end{equation*}
Thus, by (\ref{0.8Xb}), (\ref{0.9Yb}), (\ref{0.95c}) and
$1-\frac{c}{2}\in[\frac{1}{2}, 1]$,
one has $0\leqslant \mu_*\perp(d-Ex_*-Fy_*)\geqslant 0$.
This, together with (\ref{0.347}, $Ax_*+By_*-b=0$
(see (\ref{3.100b})) and Lemma \ref{Proposition1}, shows that
$(w_*, \mu_*, \alpha_*^x, \gamma_*^x, \alpha_*^y, \gamma_*^y)$
is a primal-dual KKT solution of problem (\ref{0.1}).
The proof is complete. \hfill$\Box$

Next, the iterative complexity, i.e., the maximum number of iterations of Algorithm 1 is analyzed under a given computational accuracy (error rate).
Based on the theoretical termination condition in Step 5, the following two termination criteria are considered.

{\bf Termination criterion 1} (absolute accuracy criterion):
\begin{equation}\label{error1}
  \epsilon_{k}:=\|({x}_{k+1}-x_k,  {y}_{k+1}-y_k,
  Ax_{k+1}+By_{k+1}-b)\| < \epsilon,
\end{equation}
where $\epsilon$ is an acceptable accuracy (error rate).
For ALMs, the accuracy of the feasible
approximation is usually low due to the equality constraints, i.e.,
$Ax_{k+1}+By_{k+1}$ approximates to $b$ with a low precision.
Therefore, 
under criterion  (\ref{error1}),
the number of iterations and
computational cost 
are usually high and the result
is not satisfactory for some practical problems. Therefore, similar to \cite{18,Wu2017},
based on (\ref{error1}), the following relative accuracy criterion
can be also considered for less demanding problems to save computational
cost.

{\bf Termination criterion 2} (relative accuracy criterion).
\begin{equation}\label{error2}
  \hat{\epsilon}_{k}:=\frac{\|({x}_{k+1}-x_k, {y}_{k+1}-y_k,
  Ax_{k+1}+By_{k+1}-b)\|}{\|(x_k,y_k,b)\|+1}< \epsilon.
\end{equation}

For an accumulation point $w_*=(x_*,y_*,\lambda_*)$ of $\{w_k\}$, denote
\begin{equation}\label{0.360}
  C_0= \mathcal{L}_\beta(w_0)-(f(x_*)+\theta(y_*)),
  ~\mathcal{N}_0=\min\{\rho\eta,\  \xi \}.
\end{equation}
$C_0$ denotes the difference between the initial value of the ALF and the ``optimal value''.
By Lemma \ref{lemma3.1} (ii), $C_0$ is constant and non-negative on $\Omega$.

Based on the monotonicity (\ref{0.245}), the number of iterations of
Algorithm 1 in worst-case can be obtained under the termination criteria
(\ref{error1}) and (\ref{error2}).

\begin{theorem}\label{rate5.1}
  Suppose that Assumption {\rm \ref{assumption1}} holds, and that the set $\Omega$ of accumulation points is nonempty.  
  Then for any non-negative integer $k$,
  there exists an iterative index $i_k\leqslant k$ such that
  \begin{equation}\label{dp4.15}
  \hat{\epsilon}_{i_k} \leqslant \epsilon_{i_k} \leqslant
  \sqrt{\frac{C_0}{\mathcal{N}_0}}\sqrt{\frac{1}{k+1}}.
  \end{equation}
  Therefore, for the given error rate $\epsilon$, the number
  of iterations of Algorithm 1 in worst-case is $\left[\frac{C_0}
  {\epsilon^2\mathcal{N}_0}\right]+1$, whether the termination
  criterion (\ref{error1}) or (\ref{error2}) is executed.
\end{theorem}

\textbf{Proof}\ \
Clearly, $\hat{\epsilon}_{k} \leqslant \epsilon_{k}$ holds for any
$k$, so it is sufficient to prove that the second inequality of
(\ref{dp4.15}) holds. First, from (\ref{0.245}) and (\ref{0.27}),
one has
\begin{equation*}
  \frac{ \rho\eta}{t_k} \parallel t^xd^u_k\parallel^2
  +\xi\|Ax_{k+1}+By_{k+1}-b\|^2\leqslant  \mathcal{L}_\beta(w_k)-\mathcal{L}_\beta(w_{k+1}), ~ \forall\ k\geqslant0.
\end{equation*}
This, together with $0<t_k\leqslant 1$ and $\ t_kd^u_k=({x}_{k+1}-x_k,
{y}_{k+1}-y_k)$, shows that
\begin{equation*}
  \rho\eta\|({x}_{k+1}-x_k, {y}_{k+1}-y_k)\|^2
  +\xi\|Ax_{k+1}+By_{k+1}-b\|^2\leqslant\mathcal{L}_\beta(w_k)
  -\mathcal{L}_\beta(w_{k+1}),\ \forall\ k\geqslant0.
\end{equation*}
Thus, it follows from $\mathcal{N}_0=\min\{\rho\eta, \xi\}$ and the
definition (\ref{error1}) of $\epsilon_k$ that
\begin{equation}\label{4.50A}
  \mathcal{N}_0(\epsilon_{k})^2 \leqslant
  \mathcal{L}_\beta(w_k)-\mathcal{L}_\beta(w_{k+1}),\
  \forall\ k\geqslant0.
\end{equation}
Summing the inequality (\ref{4.50A}) from $0$ to $k$, and in view of
$\mathcal{L}_\beta(w_{k+1})\geqslant \mathcal{L}_\beta(w_*)
=f(x_*)+\theta(y_*)$, one has
\begin{equation*}
\mathcal{N}_0\sum\limits_{j=0}^k (\epsilon_{j})^2 \leqslant
\mathcal{L}_\beta(w_0)-\mathcal{L}_\beta(w_{k+1}) \leqslant \mathcal{L}_\beta(w_0)-(f(x_*)+\theta(y_*))=C_0.
\end{equation*}
For any $k\geqslant0$, we consider the index $i_k$ such that
$\epsilon_{i_k}=\min\{\epsilon_0, \epsilon_1, \dots, \epsilon_k\}$.
Then the relationships above give $(k+1)\mathcal{N}_0(\epsilon_{i_k})^2
\leqslant C_0$, so (\ref{dp4.15} holds, and it further implies that
the rest conclusion of the theorem holds. $\hfill\Box$

\section{Superlinear and quadratic rates of convergence 
}
In this section, the strong convergence, superlinear and quadratic rates of convergence of Algorithm 1 are further analyzed.
As a preparation for the analysis in this section, a few basic conclusions are given below,
and the readers may refer to \cite[Corollary 1.1.8, Theorems 1.1.31 , 1.1.30]{jianjinbao2010}.
\begin{lemma}\label{lemmaN5.4}
{\rm (i)}\ If the sequence $\{v_k\}$ of real vectors has an isolated accumulation point $v_*$, and
$\lim\limits_{k\rightarrow \infty}\|v_{k+1}-v_k\|=0$, then
$\lim\limits_{k\rightarrow \infty}v_k=v_*$.

 {\rm (ii)} \ Suppose that the sequence $\{v_k\}$ of vectors converges to $v_*$,
 and has a recursive formula $v_{k+1}=v_k+d_k$, with $d_k\rightarrow 0$. If
 $\|v_{k+1}-v_*\|=o(\|v_k-v_*\|)+o(\|d_k\|)$ is satisfied, then  $\|v_{k+1}-v_*\|=o(\|v_k-v_*\|)$ holds true, i.e., $\{v_k\}$ converges superlinearly to $v_*$.

{\rm (iii)}\ If the sequence $\{v_k\}$ of vectors converges superlinearly to $v_*$,
 then $\|v_{k+1}-v_k\|\sim \|v_k-v_*\|$, i.e.,
  $\|v_{k+1}-v_k\|/\|v_k-v_*\|\rightarrow 1$.

\end{lemma}

\begin{theorem}\label{theorem2}
  Suppose that Assumptions {\rm \ref{assumption1}} and {\rm \ref{assumption2}}
  hold. If the sequence $\{w_k\}$ generated by Algorithm 1 is bounded and has an isolated accumulation point
  $w_*:=(x_*,y_*,\lambda_*)$,
   then $
  \lim\limits_{k\rightarrow \infty}d_k^u= 0,\lim\limits_{k\rightarrow \infty}w_k = w_*,\
  \lim\limits_{k\rightarrow \infty}(\alpha_k^x, \gamma_k^x, \alpha_k^y, \gamma_k^y)=
  (\alpha_*^x, \gamma_*^x, \alpha_*^y, \gamma_*^y)$ and $\lim\limits_{k\rightarrow \infty}\mu_k^x=\lim\limits_{k\rightarrow \infty} \mu_k^y=\mu_*,
$
where $(\lambda_*, \mu_*, \alpha_*^x, \gamma_*^x, \alpha_*^y, \gamma_*^y)$
is the KKT multiplier of problem (\ref{0.1}) corresponding to the
KKT point $(x_*, y_*)$. Therefore, Algorithm 1 is strong convergence.
\end{theorem}

{\bf Proof}\ \
First, notice that $\{w_k\}$ being bounded. It follows from Lemma \ref{lemma3.1} (ii)
that $\lim\limits_{k\rightarrow \infty}(Ax_k+By_k-b)=0$ and
$\lim\limits_{k\rightarrow \infty}d_k^u=0$. Then one has
\begin{equation*}\label{k+1-k}
(u_{k+1}-u_k)=
t_kd_k^u\rightarrow 0, \
(\lambda_{k+1}-\lambda_k)=
\xi(Ax_{k+1}+By_{k+1}-b)\rightarrow 0, {\rm as}~ k\rightarrow \infty.
\end{equation*}
Therefore, $\|w_{k+1}\!-\!w_k\|\rightarrow 0$.
This, along with the assumption that  $w_*$ being an isolated accumulation point of
$\{w_k\}$ and Lemma \ref{lemmaN5.4} (i), shows that
$\lim\limits_{k\rightarrow \infty}w_k \! = \!w_*$.
Second, by Lemma \ref{lemma3.2}, it is known that  $\{(\mu_k^x, \mu_k^y, \alpha_k^x,
\gamma_k^x,$ $\alpha_k^y, \gamma_k^y)\}$ is bounded.
From Theorem \ref{theorem1}, for each accumulation point
$(\tilde{\mu}_*^x, \tilde{\mu}_*^y, \tilde{\alpha}_*^x,
\tilde{\gamma}_*^x, \tilde\alpha_*^y, \tilde\gamma_*^y)$ of
$\{(\mu_k^x, \mu_k^y, \alpha_k^x, \gamma_k^x,\alpha_k^y, \gamma_k^y)\}$,
one has $\tilde{\mu}_*^x=\tilde{\mu}_*^y=:\tilde{\mu}_*$.
Further, $(\tilde{\mu}_*, \tilde{\alpha}_*^x, \tilde{\gamma}_*^x,
\tilde\alpha_*^y, \tilde\gamma_*^y)$ together with $w_*$ constitutes a
primal-dual KKT solution of problem (\ref{0.1}.
Also by Assumption \ref{assumption2}, one knows that the KKT multiplier
corresponding to the accumulation point $w_*$ is unique.
Therefore, the multiplier sequence $\{(\mu_k^x, \mu_k^y, \alpha_k^x,
\gamma_k^x, \alpha_k^y, \gamma_k^y)\}$ has a unique accumulation point
$(\mu_*, \alpha_*^x, \gamma_*^x, \alpha_*^y, \gamma_*^y)$, and thus
converges to this unique accumulation point. \hfill$\Box$

A second-order sufficient condition for $\{w_k\}$ to has an isolated
accumulation point is given below, which also plays an important
role in the analysis of superlinear convergence.

\begin{assumption}\label{assumption3}
  {\bf (i)} The objective functions $f(x)$ and $\theta(y)$ are at
  least second-order continuously differentiable in a neighborhood
  of the partially feasible set $\mathcal{F}_+$.

  {\bf (ii)} The sequence $\{w_k\}$ generated by Algorithm 1 is
  bounded, and its accumulation point $w_*=(u_*,\lambda_*)$ satisfies the following
  two requirements.

 \ \  {\bf (iia)}
 The linearly independent constraint qualification {\rm (LICQ)} is satisfied at $ u_*=(x_*,y_*)$,
 i.e., the corresponding gradient matrix of active constraint
  \begin{equation*}
    \mathcal{A}_*:=\left(
    \begin{array}{cc}
      A&B\\
    E_{S_{u_*}}&F_{S_{u_*}}\\
    C_{I_{x_*}} &0\\
    0&D_{J_{y_*}}
    \end{array}
    \right)=\left(
      \begin{array}{cc}
        A~~~~~~B\\
        ~~\mathcal{B}_{(x_*,y_*)}
      \end{array}
      \right)
    \end{equation*}
is full row rank.

  \ \ {\bf (iib)}
By Theorem {\rm \ref{theorem1}} and the LICQ
  above, there exists a unique multiplier
  $(\mu_*, \alpha^x_*, \gamma^x_*, \alpha^y_*, \gamma^y_*)$ such that
  $w_*$ together with this 
  constitutes a primal-dual KKT solution of problem (\ref{0.1}).
  Suppose that the strong second-order sufficient (SSOS) condition at
  $(w_*, \mu_*, \alpha^x_*, \gamma^x_*, \alpha^y_*, \gamma^y_*)$ is satisfied, i.e.,
 \begin{equation}\label{SSOSC1}
  {d^u}^\top
  \nabla^2 F(x_*,y_*) d^u>0,
  \forall~ 0\neq d^u \in M_*:=\{d^u 
  \in \mathbb{R}^{n_1+n_2}:\ \mathcal{A}^+_*d^u=0\},
  \end{equation}
where,
  \begin{equation*}
    \mathcal{A}^+_*=\left(
    \begin{array}{cc}
      A&B\\
    E_{S_*^+}&F_{S_*^+}\\
    C_{I_*^+} &0\\
    0&D_{J_*^+}
    \end{array}
    \right),
    \ \  \mathcal{B}^+_*=\left(
    \begin{array}{cc}
    E_{S_*^+}&F_{S_*^+}\\
    C_{I_*^+} &0\\
    0&D_{J_*^+}
    \end{array}
    \right),
  \end{equation*}
 and
  \vspace*{-0.3cm}
  \begin{equation}\label{0.349}
  \begin{array}{ll}
  I_*^{l_+}= \{i\in I_*^l:=I_{x_*}^l:(\alpha_*^x)_i> 0\},
  ~I_*^{v_+}= \{i\in I_*^v:=I_{x_*}^v:(\gamma_*^x)_i> 0\},\\
  J_*^{s_+}= \{j\in J_*^s:=J_{y_*}^s:(\alpha_*^y)_j> 0\},
  ~J_*^{r_+}= \{j\in J_*^r:=J_{y_*}^r:(\gamma_*^y)_j> 0\},\\
  I_*^+=I_*^{l_+}\cup I_*^{v_+},\ J_*^+=J_*^{s_+}\cup J_*^{r_+},\ S_*^+=\{i\in S_*:=S_{u_*}:(\mu_*)_i>0\}.
  \end{array}
  \end{equation}
\end{assumption}

\begin{theorem}\label{guli}
Suppose that Assumptions {\rm \ref{assumption1}} and {\rm \ref{assumption3}} hold.
Then the sequence $\{w_k\}$ generated by Algorithm 1 has an isolated accumulation point
$w_*$, and $\lim\limits_{k\rightarrow \infty}(w_k,\mu_k^x,\mu_k^y,
\alpha^x_k, \gamma^x_k, \alpha^y_k, \gamma^y_k)
=(w_*,\mu_*,\mu_*,\alpha^x_*, \gamma^x_*, \ \alpha^y_*, \gamma^y_*)$.
\end{theorem}

{\bf Proof}\ \
First, by Assumption \ref{assumption3}, it follows from
\cite[Corollary 1.4.3]{jianjinbao2010} that $(x_*,y_*)$ is an isolated
KKT point of problem (\ref{0.1}), and the corresponding multiplier
$(\lambda_*, \mu_*, \alpha^x_*, \gamma^ x_*, \alpha^y_*, \gamma^y_*)$
is unique. Furthermore, it is not difficult to know that
$w_*=(x_*, y_*, \lambda _*)$ is an isolated accumulation point of $\{w_k\}$.
Otherwise, there exists an infinite sequence
$\{\hat{w}_i=(\hat{x}_i,\hat{y}_i,\hat{\lambda}_i)\}$ of accumulation points of
$\{w_k\}$ such that $\hat{w}_i \rightarrow w_*$ and $\hat{w}_i \neq w_*$.
Thus, by Theorem \ref{theorem1}, $\{(\hat{x}_i,\hat{y}_i)\}$ is a sequence
of KKT points of problem (\ref{0.1}, and $(\hat{x}_i,\hat{y}_i)\not=
(x_*,y_*)$ by the uniqueness of the corresponding multiplier.
This contradicts the fact that $(x_*,y_*)$ is an isolated KKT point
of problem (\ref{0.1}).  Finally, in view of the isolation of the accumulation point
$w_*$ and the fact that Assumption \ref{assumption3} (ii) implies Assumption \ref{assumption2}, by Theorem {\rm \ref{theorem2}}, one has
$\lim\limits_{k\rightarrow \infty}(w_k,\mu_k^x,\mu_k^y, \alpha^x_k, \gamma^x_k, \alpha^y_k, \gamma^y_k)
=(w_*,\mu_*,\mu_*,\alpha^x_*, \gamma^x_*, \alpha^y_*, \gamma^y_*)$. $\hfill\Box$

\begin{lemma}\label{lem5.1}
  {\rm Under the SSOS condition
  (\ref{SSOSC1},  there exists a constant  $\beta_*>0$ such that
$\nabla^2_{uu}\mathcal{L}_{\beta}(w^*)$ is positive definite on the
null space $M_*^+:=\{d^u \in \mathbb{R}^{n_1+n_2}:\
\mathcal{B}^+_*d^u=0\}$ for all $\beta\geqslant\beta_*$.}
\end{lemma}

{\bf{Proof}}~~Suppose by contradiction that the claim is not true.
Then there exists an infinite sequence $\{\beta_k\}\rightarrow +\infty$
and associated sequence $\{d^u_k\not=0\}\subseteq M^+_*$ such that
${d^u_k}^\top\nabla^2_{uu} \mathcal{L}_{\beta_k}(w^*)d^u_k\leqslant 0$.
It follows from the definition (\ref{La}) of $\mathcal{L}_{\beta}(w)$ that
\begin{equation}\label{Luub}
        \nabla_{uu}^2\mathcal{L}_\beta(w)
    =\nabla^2F(x,y)+\beta (A\  B)^\top (A\  B).
\end{equation}
Therefore, one has
    \[
    0\geqslant{d^u_k}^\top\nabla^2_{uu} \mathcal{L}_{\beta_k}(w^*)d^u_k={d^u_k}^\top\nabla^2F(x_*,y_*)d^u_k+\beta_k\|(A\  B)d^u_k\|^2.
    \]
Denote $\bar{d}^u_k=\frac{d^u_k}{\|d^u_k\|}$. Then, without loss of generality,
it can be assumed that $\bar{d}^u_k\rightarrow \bar{d}^u$ as $k\rightarrow\infty$.
Obviously, $0\not=\bar{d}^u\in M^+_*$.
Therefore, it follows from the relationship above that
    \[
    b_k:=(\bar{d}^u_k)^\top\nabla^2F(x_*,y_*)\bar{d}^u_k+\beta_k\|(A\ \ B) \bar{d}^u_k\|^2\leqslant 0,\ \forall \
    k.
    \]
Dividing the above inequality by $\beta_k$ and taking the limit of
$k\rightarrow\infty$, one has $(A~~B)\bar{d}^u=0$. This, together
with $\bar{d}^u\in M^+_*$ shows that $\bar{d}^u\in M_*$ (see (\ref{SSOSC1})).
Therefore, in view of the positive definiteness of $\nabla^2F(x_*,y_*)$  on $M_*$, one has 
\[
 b_k\geqslant(\bar{d}^u_k)^\top\nabla^2F(x_*,y_*)\bar{d}^u_k\rightarrow (\bar{d}^u)^\top\nabla^2F(x_*,y_*)\bar{d}^u>0.
\]
This contradicts $b_k\leqslant 0$, and the proof is complete. \hfill$\Box$

For convenience of presentation, denote Algorithm 1 corresponding to $c=c_0$ as the
PDF-SQO$_{c_0}$ method. From Steps 2-3 of the PDF-SQO$_{1}$ method,
we know that $c_{\rm max}=1$, and the initial test value in
Step 4 reaches $1$. Moreover, we will show that the PDF-SQO$_{1}$ method
can achieve a step size $1$ when the equality constraint (\ref{0.1b})
satisfies $A^\top B=0$, i.e., it can overcome the
Maratos effect \cite{Maratos}. If $A^\top B=0$, from (\ref{**}),
(\ref{Hku}) and (\ref{La}), it follows that
  \begin{equation}\label{Luua}
  \tilde{\mathcal{H}}^u_k=\hat{\mathcal{H}}^u_k
  =\mathcal{H}^u_k ={\rm diag}(\mathcal{H}_k^x, \mathcal{H}_k^y)
  ={\rm diag}(H_k^x+\beta A^\top A, H_k^y+\beta B^\top B),
  \end{equation}
  \begin{equation}\label{Luub}
        \nabla_{uu}^2\mathcal{L}_\beta(w_k)
        =\left(
          \begin{array}{cc}
          \nabla^2f(x_k)+\beta A^\top A & \\
           & \nabla^2\theta(y_k)+\beta B^\top B
          \end{array}
          \right).
  \end{equation}
Denote active sets by
  \begin{equation}\label{A-sets}
  \begin{array}{ll}
    I^l_k=I^l_{\tilde{x}_{k+1}},\ I^v_k=I^v_{\tilde{x}_{k+1}},\ J^s_k=J^s_{\tilde{y}_{k+1}},\ J^r_k=J^r_{\tilde{y}_{k+1}},\\ I_k=I_{\tilde{x}_{k+1}},\ J_k=J_{\tilde{y}_{k+1}},\ S_k=S_{(\tilde{x}_{k+1},\tilde{y}_{k+1})},
    \end{array}
  \end{equation}
  \begin{equation}\label{hatS^x_k}
  \begin{array}{ll}
     \hat{S}^x_k=\{i\in S:\ (E\tilde{x}_{k+1}+Fy_k-d-\frac12 h_k)_i=0\},\\
     \hat{S}^y_k=\{i\in S:\ (Ex_k+F\tilde{y}_{k+1}-d-\frac12 h_k)_i=0\},
     \end{array}
  \end{equation}
  and
  \begin{equation}\label{S^+_k}
    S_k^+=
    \left\{
      \begin{array}{ll}
    S_k, & {\rm if}\ \tilde{u}_{k+1}\
    {\rm is\ generated\ by\ Step\ 3,\ i.e.,
    \  QO \ (\ref{0.35}) }; \\
    \hat{S}^x_k\cap \hat{S}^y_k,  & {\rm if}\ \tilde{u}_{k+1}
    \ {\rm is\ generated\ by\ Step\ 1,\ i.e.,
    \  QOs}\ (\ref{0.6X}\ {\rm and}\ (\ref{0.7Y}).
      \end{array}
    \right.
  \end{equation}
  Define matrices
  \begin{equation}\label{Hk^}
    P_k= E_{n_1+n_2}-\mathcal{B}_k^\top(\mathcal{B}_k\mathcal{B}_k^\top)^{-1}\mathcal{B}_k,
    ~\mathcal{B}_k=
    \left(
      \begin{array}{cc}
      E_{S^+_k} &  F_{S^+_k}\\
      C_{I_k} & 0\\
      0 & D_{J_k}
    \end{array}
    \right).
  \end{equation}
By relationship (\ref{setseqc}) later and Assumption
\ref{assumption3} (iia), the matrix $\mathcal{B}_k$ defined above
is full row rank (when $k$ is sufficiently large), so the above
projection operation $_k$ is well-defined.

First, to ensure that the unit step is accepted by the
PDF-SQO$_{1}$ method, the matrix $\mathcal{H}^u_k$ needs to
satisfy the following assumption.
\begin{assumption}\label{assumption4}
Suppose that the equality constraint (\ref{0.1b}) satisfies
$A^\top B=0$ $($this automatically holds when the equality
constraint vanishes), 
and assume that the sequence $\{\mathcal{H}^u_k\}$ of matrices
generated by PDF-SQO$_{1}$ 1 satisfies one of
the following two conditions:
\begin{equation}\label{tiaojian1}
\|(\mathcal{H}^u_k-\nabla_{uu}^2\mathcal{L}_\beta(w_k))d^u_k\|=o(\|d^u_k\|);
\end{equation}
\vspace*{-0.2cm}
\begin{subequations}\label{tiaojian2}
 \begin{numcases}{}
          \|P_k(\mathcal{H}^u_k-\nabla_{uu}^2\mathcal{L}_\beta(w_k))d^u_k\|=o(\|d^u_k\|),\label{tiaojian2a}\\
          (\alpha^x_*)_{I_*^l}>0, (\gamma^x_*)_{I_*^v}>0, (\alpha^y_*)_{J_*^s}>0, (\gamma^y_*)_{J_*^r}>0, (\mu_*)_{S_*}>0.\label{tiaojian2b}
 \end{numcases}
\end{subequations}
\end{assumption}

By (\ref{**}), (\ref{Hku}) and (\ref{Luub}) as well as $A^\top B=0$,
a strong form of (\ref{tiaojian1}) and (\ref{tiaojian2a}) is
\begin{equation}\label{tiaojian3}
    \lim_{k\rightarrow\infty}(H_k^x-\nabla^2f(x_k))=0,\
    \lim_{k\rightarrow\infty}(H_k^y-\nabla^2\theta(y_k))=0.
\end{equation}
This is consistent with Step 5 of the proposed method.
\begin{lemma}\label{lemma5}
Suppose that Assumptions {\rm \ref{assumption1}} and
{\rm \ref{assumption3}} hold. Then, when $k$ is sufficiently large,

{\rm (i)} the active sets generated by the {\rm PDF-SQO$_{1}$}
method satisfy the following relationships:
\begin{subequations}\label{setseq}
 \begin{numcases}{}
 I^{l+}_*\subseteq I_k^l \subseteq I^l_*, \ I^{v+}_*\subseteq
 I_k^v \subseteq I^v_*,\ J^{s+}_*\subseteq J_k^s \subseteq
 J^s_*,\ J^{r+}_*\subseteq J_k^r \subseteq J^r_*, \label{setseqa} \\
 S^+_*\subseteq \hat{S}^x_k \subseteq S_*,\ S^+_*\subseteq
 \hat{S}^y_k \subseteq S_*,\ S^+_*\subseteq S_k \subseteq S_*,\label{setseqb}\\
 I^+_*\subseteq I_k \subseteq I_*:=I_{x^*}, \ J^+_*\subseteq
 J_k \subseteq J_*:=J_{y^*},\  S^+_*\subseteq S_k^+ \subseteq
 S_*; \label{setseqc}
 \end{numcases}
\end{subequations}

{\rm (ii)} if the strict complementarity condition in
(\ref{tiaojian2b}) also holds, then the relation $``\subseteq"$
in (\ref{setseq}) can all be strengthened to $``="$.
\end{lemma}

{\bf Proof} \ \
(i) First, for $i\in I_*^{l+}$, one has $(\alpha_k^x)_i
\rightarrow (\alpha_*^x)_i>0$ by Theorem \ref{guli}.
Further, $(\alpha_k^x)_i>0$ holds when $k$ is sufficiently large,
thus $i\in I^l_k$. So $I_*^{l+}\subseteq I_k^l$ holds.
Second, one can easily verify that $I_k^l\subseteq I_*^l$ holds for sufficiently large $k$.
If not, there exist an infinite set $K$ and a fixed $i$ such that
$i\in I_k^l\backslash I_*^l, \forall~ k\in K$.
Thus, $l_i=c_i\tilde{x}_{k+1}=c_i(x_k+d_k^x) \rightarrow
c_ix_*, i\in I_*^l$, which is a contradiction.
Therefore, the relations $I^{l+}_*\subseteq I_k^l \subseteq
I^l_*$ are proved. The rest of the proofs of the relations in
(\ref{setseqa} can be proved similarly.

For $i\in S_*^{+}$, one has $(\mu_k^x)_i \rightarrow (\mu_*)_i>0$
by Theorem \ref{guli}. Thus, $(\mu_k^x)_i>0$ holds when $k$ is
sufficiently large. Further, by the complementary condition
(\ref{0.8Xb}), one knows that $i\in \hat{S}^x_k$.
Therefore $S_*^{+}\subseteq \hat{S}^x_k$ holds.
Second, we prove that $\hat{S}^x_k \subseteq S_*$ holds for
sufficiently large $k$. If not, there exist an infinite set $K$
and a fixed $j$ such that $j\in \hat{S}^x_k\backslash S_*,
\forall~ k\in K$. So $(Ex_*+Fy_*-d)_j=\lim_{k\in K}(E\tilde{x}_{k+1}
+Fy_k-d-\frac12 h_k)_j=0$, which means $j\in S_*$, leading to a
contradiction. Therefore, the relations $S^+_*\subseteq
\hat{S}^x_k \subseteq S_*$ are proved. The rest of the proofs of
the relations in (\ref{setseqb}) can be proved similarly, and
(\ref{setseqc}) is a direct corollary of the conclusions
(\ref{setseqa}) and (\ref{setseqb}).

(ii) If the strict complementarity condition in (\ref{tiaojian2b})
holds, then $I^{l+}_*=I^l_*, \ I^{v+}_*= I^v_*,\ J^{s+}_*= J^s_*,
\ J^{r+}_*= J^r_*,\ I^+_*=I_* \ J^+_*= J_*,\ S^+_*= S_*.$
Thus, the relations in (\ref{setseq}) show that the conclusion holds. \hfill$\Box$

\begin{theorem}\label{theorem3}
 Suppose that Assumptions {\rm \ref{assumption1}}, {\rm \ref{assumption3}}
 and {\rm \ref{assumption4}} hold. If the parameter $\rho\in (0,0.5)$,
 then, when the iteration index $k$ is sufficiently large,
 the step size of PDF-SQO$_{1}$ method reaches $1$, i.e., $t_k\equiv1$.
\end{theorem}

{\bf Proof}\ \
First, $c_{\max}=1$ holds from Steps 2-3. Thus, by Step 4,
it is sufficient to prove that the inequality (\ref{0.24})
holds for $t=1$ and all sufficiently large $k$. Further,
from the second-order Taylor expansion, one has
  \begin{equation}\label{zhixu}
    \begin{array}{ll}
      \nu_k&:=\mathcal{L}_\beta(u_k+d_k^u,\lambda_k)-\mathcal{L}_\beta(u_k,\lambda_k)+\rho\|d_k^u\|_{\tilde{\mathcal{H}}_k^u}^2\\
      &=\nabla_u \mathcal{L}_\beta(w_k)^\top d_k^u+\frac{1}{2}{d^u_k}^\top \nabla_{uu}^2\mathcal{L}_\beta(w_k)d^u_k+\rho\|d^u_k\|_{\tilde{\mathcal{H}}^u_k}^2+o(\|d_k^u\|^2).
    \end{array}
  \end{equation}

Next, it is further proved that $\nu_k\leqslant0$ holds under
the conditions (\ref{tiaojian1}) and (\ref{tiaojian2}), respectively.

(i) Suppose that (\ref{tiaojian1}) holds. Then, by (\ref{0.20}), (\ref{zhixu}),
$\tilde{\mathcal{H}}^u_k=\mathcal{H}_k^u$, (\ref{tiaojian1}) and
(\ref{yizhi}), one has
\begin{equation*}\label{nuk}
 \begin{array}{ll}
 \nu_k&\leqslant -\|d_k^u\|^2_{\mathcal{H}_k^u}
 +\frac{1}{2}{d_k^u}^\top \nabla^2_{uu}\mathcal{L}_\beta(w_k)d_k^u
 +\rho\|d_k^u\|^2_{\mathcal{H}_k^u}+o(\|d_k^u\|^2)\\
 &=(\rho-\frac{1}{2})\|d_k^u\|^2_{\mathcal{H}_k^u}
 +\frac{1}{2}{d_k^u}^\top (\nabla^2_{uu}\mathcal{L}_\beta(w_k)-\mathcal{H}_k^u)d_k^u
 +o(\|d_k^u\|^2)\\
 &\leqslant\eta (\rho-\frac{1}{2})\|d^u_k\|^2+
 o(\|d_k^u\|^2)\leqslant 0.
 \end{array}
\end{equation*}

(ii) Suppose that (\ref{tiaojian2}) holds. First, from $c=1$,
(\ref{xiajiang}), (\ref{0.951}) (notice that $\mu_k^x=\mu_k^y$
in Step 3), $\hat{\mathcal{H}}_k^u=\mathcal{H}_k^u$ and Lemma
\ref{lemma5} (ii), it is easy to see that
\begin{equation*}\label{leq*}
 \begin{array}{ll}
 \nabla_u\mathcal{L}_\beta(w_k)^\top d_k^u&\leqslant\!-
 (Cx_k\!-\!l)_{I_k^l}^\top(\alpha_k^x)_{I_k^l}\!-\!
 (v\!-\!Cx_k)_{I_k^v}^\top(\gamma_k^x)_{I_k^v}\!-\!
 (Dy_k-s)_{J_k^s}^\top(\alpha_k^y)_{J_k^s}\!-\\
 &~~~~~(r\!-\!Dy_k)_{J_k^r}^\top(\gamma_k^y)_{J_k^r}+
 \frac{1}{2}(h_k^\top)_{S^+_k}(\mu_k^x+\mu_k^y)_{S^+_k}\!-\!
 \|d_k^u\|_{{\mathcal{H}}_k^u}^2.
 \end{array}
\end{equation*}
This, along with (\ref{zhixu}) and $\tilde{\mathcal{H}}_k^u=\mathcal{H}_k^u$, shows that   \begin{equation*}\label{nuk0}
 \begin{array}{ll}
    \nu_k&\leqslant (\rho-\frac{1}{2})\|d_k^u\|^2_{\mathcal{H}_k^u}
    +\frac{1}{2}{d_k^u}^\top (\nabla^2_{uu}\mathcal{L}_\beta(w_k)-\mathcal{H}_k^u)d_k^u
    -(Cx_k-l)_{I_k^l}^\top(\alpha_k^x)_{I_k^l}-\\
    &\ \ \ \ \ \ (v-Cx_k)_{I_k^v}^\top(\gamma_k^x)_{I_k^v}-
    (Dy_k-s)_{J_k^s}^\top(\alpha_k^y)_{J_k^s}-
    (r-Dy_k)_{J_k^r}^\top(\gamma_k^y)_{J_k^r}+\\
    &\ \ \ \ \ \ \frac{1}{2}(h_k^\top)_{S^+_k}(\mu_k^x+\mu_k^y)_{S^+_k}
    +o(\|d_k^u\|^2).
 \end{array}
\end{equation*}
Let
\begin{equation*}
 \sigma_*=\frac{1}{2}\min\{(\alpha_*^x)_i,i\in I_*^l;
 (\gamma_*^x)_i, i\in I_*^v;(\alpha_*^y)_j, j\in J_*^s;
 (\gamma_*^y)_j, j\in J_*^r; (\mu_*)_i, i\in S_*\} \overset{(\ref{tiaojian2b}}{>}0.
\end{equation*}
Then, the two relationships above, together with $(\mu_k^x,\mu_k^y, \alpha^x_k, \gamma^x_k, \alpha^y_k, \gamma^y_k)\rightarrow
  (\mu_*,\mu_*,\alpha^x_*, \gamma^x_*, \alpha^y_*, \gamma^y_*)$ and (\ref{setseqa})-(\ref{setseqb}), provide
\begin{equation}\label{vk3}
  \begin{array}{ll}
\nu_k \leqslant (\rho-\frac{1}{2})\|d_k^u\|^2_{\mathcal{H}_k^u}
+\frac{1}{2}{d_k^u}^\top (\nabla^2_{uu}\mathcal{L}_\beta(w_k)-\mathcal{H}_k^u)d_k^u
-\sigma_*\|z_k\|_1+o(\|d_k^u\|^2),
\end{array}
\end{equation}
where
\begin{equation}\label{z_k}
  z_k^\top=(
  (-h_k)_{S^+_k}^\top,\
  (l-Cx_k)_{I_k^l}^\top, \
  (v-Cx_k)_{I_k^v}^\top,\
  (s-Dy_k)_{J_k^s}^\top, \
  (r-Dy_k)_{J_k^r}^\top).
\end{equation}

On the other hand, by the definitions (\ref{hatS^x_k}) and (\ref{S^+_k}) of the active sets, for each $i\in S^+_k$, one has
\begin{equation}\label{N01}
  \begin{array}{ll}
  e_id^x_k+f_id^y_k & =(e_i\tilde{x}_{k+1}+f_i\tilde{y}_{k+1})-
  (e_ix_k+f_iy_k)\\
  & =
  \left\{
  \begin{array}{l}
  d_i-(e_ix_k+f_iy_k)=(-h_k)_i, \ {\rm if}\ S^+_k=S_k; \\
  (\frac12 (h_k)_i+d_i-f_iy_k)+(\frac12 (h_k)_i+d_i-e_ix_k)-
  (e_ix_k+f_iy_k)\\
  =(-h_k)_i, \ {\rm if}\ S^+_k=\hat{S}^x_k\cap\hat{S}^y_k.
  \end{array}
  \right.
  \end{array}
\end{equation}
Similarly, one also has
\begin{equation}\label{N05}
  \left\{
  \begin{array}{l}
  C_{I_k}d^x_k=C_{I_k}\tilde{x}_{k+1}-C_{I_k}x_k= ((l-Cx_k)_{I_k^l}^\top, \
  (v-Cx_k)_{I_k^v}^\top)^\top,\\
  D_{J_k}d^y_k=D_{J_k}\tilde{y}_{k+1}-D_{J_k}y_k=((s-Dy_k)_{J_k^s}^\top, \ (r-Dy_k)_{J_k^r}^\top)^\top.
  \end{array}
  \right.
\end{equation}
Therefore, from the relationships (\ref{Hk^}) and (\ref{z_k})-(\ref{N05}), it follows that
\ $\mathcal{B}_kd_k^u=z_k$. This, together with (\ref{Hk^}), further implies that
\begin{equation}\label{d^uk}
  d_k^u
  =P_kd_k^u+\mathcal{B}_k^\top(\mathcal{B}_k\mathcal{B}_k^\top)^{-1}
  \mathcal{B}_kd_k^u, \ \ \mathcal{B}_kd_k^u=z_k.
\end{equation}
Furthermore, from (\ref{N01}) and (\ref{N05}), one has
\begin{equation}\label{B_kx}
  E_{S^+_k}\tilde{x}_{k+1}+F_{S^+_k}\tilde{y}_{k+1}=d_{S^+_k},\
  C_{I_k}\tilde{x}_{k+1}
  =\left(
  \begin{array}{c} l_{{I_k^l}}\\
  v_{{I_k^v}}
  \end{array}
  \right),\ \
  D_{J_k}\tilde{y}_{k+1}
  =\left(
  \begin{array}{c} r_{{J_k^r}}\\
  s_{{J_k^s}}
  \end{array}
  \right).
\end{equation}
Thus, it follows from (\ref{d^uk}) and (\ref{tiaojian2a}) that
\begin{equation*}
     \begin{array}{ll}
      {d_k^u}^\top (\nabla^2_{uu}\mathcal{L}_\beta(w_k)-\mathcal{H}_k^u)d_k^u
      &\ \overset{(\ref{d^uk}}{=} {d^u_k}^\top(\nabla_{uu}^2\mathcal{L}_\beta(w_k)
      -\mathcal{H}^u_k)
        (P_kd^u_k+\mathcal{B}_k^\top(\mathcal{B}_k
        \mathcal{B}_k^\top)^{-1}
  z_k)\\
      &\ \ \ ={d^u_k}^\top(\nabla_{uu}^2\mathcal{L}_\beta(w_k)
      -\mathcal{H}^u_k)P_kd^u_k+O(\|d^u_k\|\cdot\|z_k\|_1)\\
      &\ \overset{(\ref{tiaojian2a}}{=}o(\|d^u_k\|^2)+o(\|z_k\|_1).
     \end{array}
\end{equation*}
Substituting the above relation into (\ref{vk3}), and combining
(\ref{0.27}) and $\rho\in(0, 0.5)$, one has
\begin{equation*}
    \begin{array}{ll}
  \nu_k\leqslant(\rho-\frac{1}{2})\eta\|d^u_k\|^2+o(\|d^u_k\|^2)-\sigma_*\|z_k\|_1+o(\|z_k\|)
  \leqslant0.
  \end{array}
\end{equation*}
The proof of the theorem is complete. \hfill$\Box$

Based on relationship (\ref{setseqc}), by \cite[Theorem 1.1.10 (2)]
{jianjinbao2010}, the following lemma holds true.
\begin{lemma}\label{lemma6}
Suppose that Assumptions {\rm \ref{assumption1}} and {\rm \ref{assumption3}} hold.
If the penalty parameter $\beta\geqslant\beta_*$ which is defined in Lemma \ref{lem5.1},
then the matrix
\begin{equation*}
\mathcal{G}_k^{*\beta}:=\left(
\begin{array}{cc}
P_k\nabla^2_{uu}\mathcal{L}_\beta(w_*)&\mathcal{B}_k^\top \\
\mathcal{B}_k &0
\end{array}
\right)
\end{equation*}
is uniform invertible, i.e., there exists a constant $\mathcal{M}>0$ such that
$\|(\mathcal{G}_k^{*\beta})^{-1}\|\leqslant \mathcal{M}$ holds for all
sufficiently large $k$. Specially, if the equality constraint
$Ax+By-b=0$ in problem (\ref{0.1}) disappears, then the claim above
holds independent with $\beta$ for matrix
\begin{equation*}
\mathcal{G}_k^{*}:=\left(
\begin{array}{cc}
P_k\nabla^2F(u_*)&\mathcal{B}_k^\top \\
\mathcal{B}_k &0
\end{array}
\right).
\end{equation*}
\end{lemma}

As the core conclusion of this section, the PDF-SQO$_{1}$ method is discussed
to achieve superlinear rate of convergence for solving a special class of problem (\ref{0.1}),
i.e., the equality constraint $Ax+By-b=0$ disappears. Denote the sub-model of (\ref{0.1}) without the
equality constraint $Ax+By-b=0$ as problem (\ref{0.1})-eq-free.
According to Lemma \ref{lemma2}, the PDF-SQO$_{c}$ method is a
completely feasible (decreasing) distributed SQO method (denoted as
{\rm CFD-SQO$_c$-eq-free} method) when it is applied to problem
(\ref{0.1})-eq-free.
\begin{theorem}\label{theorem4}
Suppose that Assumptions {\rm \ref{assumption1}}, {\rm \ref{assumption3}}
and {\rm \ref{assumption4}} hold, and parameters $\rho\in(0,0.5)$ and $\tau_1>1$.
Then the sequence $\{u_k=(x_k,y_k)\}$
generated by the {\rm CFD-SQO$_1$-eq-free} method superlinearly converges
to $u^*=(x_*,y_*)$, i.e., $\|u_{k+1}-u_*\|=o(\|u_k-u_*\|)$.
Namely, the {\rm CFD-SQO$_1$-eq-free} method is superlinearly convergent.
\end{theorem}
{\bf Proof}\ \
First, in view of the fact that the objective function $F(u)=f(x)+g(y)$
is at least second-order continuously differentiable,
taking the gradient map $\nabla F(u)=(\nabla f(x)^\top,
\nabla \theta(y)^\top)^\top$  of $F(u)$ at $u_*$ and Taylor
expansion, one has
\begin{equation}\label{F1}
 \nabla F(u_k)=\nabla F(u_*)+\nabla^2 F(u_*)(u_k-u_*)+\varphi_k,
\end{equation}
where $\varphi_k$ is given by
\begin{equation}\label{F1b}
\varphi_k=o(\|u_k-u_*\|)\ \ {\rm if}\ F(u)\in {\rm C}^2; \ \varphi_k=O(\|u_k-u_*\|^2)\ \ {\rm if}\ F(u)\in {\rm C}^3.
\end{equation}
Again, in view of the equality constraint $Ax+By-b=0$ disappearing
 and relationship (\ref{setseqc}), and taking into account
$(u^*,\mu_*, \alpha_*^x, \gamma_*^x, \alpha_*^y, \gamma_*^y)$
satisfying the KKT condition (\ref{0.17}), one has
\begin{equation}\label{F2}
  \nabla F(u_*)=
  \left(
    \begin{array}{c}
      \nabla f(x_*) \\
      \nabla \theta(y_*)
    \end{array}
  \right)=-\mathcal{B}_k^\top
  \left(
    \begin{array}{c}
      (\mu_*)_{S^+_k}\\
      (\gamma_*^x-\alpha_*^x)_{I_k}\\
      (\gamma_*^y-\alpha_*^y)_{J_k}
    \end{array}
  \right).
\end{equation}
The two equalities (\ref{F1}) and (\ref{F2}) above, together with $P_k\mathcal{B}_k^\top=0$, show that
\begin{equation}\label{P_kF_k}
\begin{array}{ll} P_k\nabla F(u_k)&
=P_k\nabla F(u_*)+P_k\nabla^2 F(u_*)(u_k-u_*)+O(\varphi_k)\\
&=P_k\nabla^2 F(u_*)(u_k-u_*)+O(\varphi_k).
\end{array}
\end{equation}

On the other hand, in view of $\tilde{\mathcal{H}}_k^u=\mathcal{H}_k^u$
and the equality constraint $Ax+By-b=0$ vanishing,
we can obtain the following relationship from (\ref{0.28a}),
\begin{equation}\label{0.28adengjia}
  \left(\!
  \begin{array}{c}
    \nabla f(x_k) \\
    \nabla \theta(y_k)
  \end{array}\!
  \right)\!
  +\!\mathcal{H}_k^ud_k^u\!+\!
  \left(
    \begin{array}{c}
      0\\
      F^\top(\mu_k^y-\mu_k^x)
    \end{array}
  \right)\!
  =-\mathcal{B}_k^\top
  \left(\!
  \begin{array}{c}
      (\mu_k^x)_{S^+_k}\\
      (\gamma_k^x-\alpha_k^x)_{I_k}\\
      (\gamma_k^y-\alpha_k^y)_{J_k}
  \end{array}\!
  \right).
\end{equation}
Therefore, combining $P_k\mathcal{B}_k^\top=0$, relations (\ref{0.28adengjia}),
(\ref{P_kF_k}) and (\ref{0.28f}), it follows that
\begin{equation*}
  \begin{array}{ll}
    P_k\mathcal{H}_k^ud_k^u&
    \overset{(\ref{0.28adengjia}}{=}-P_k\nabla F(u_k)-P_k
    \left(
      \begin{array}{c}
        0 \\
        F^\top(\mu_k^y-\mu_k^x)
      \end{array}
      \right)\\
      & \overset{(\ref{P_kF_k}}{=}
      -P_k\nabla^2F(u_*)(u_k-u_*)-P_k
      \left(
        \begin{array}{c}
          0 \\
          F^\top(\mu_k^y-\mu_k^x)
        \end{array}
      \right)
      +O(\varphi_k)\\
      &\overset{(\ref{0.28f}}{=}-P_k\nabla^2F(u_*)(u_k-u_*)
      +O(\varphi_k)+O(\|d_k^u\|^{\tau_1}).
  \end{array}
\end{equation*}
This, along with $ u_{k+1}=u_k+d_k^u$, further yields
\begin{equation*}\label{Qk}
  \begin{array}{ll}
    P_k\nabla^2F(u_*)(u_{k+1}-u_*)
    &=P_k(\nabla^2F(u_*)-\mathcal{H}_k^u)d_k^u
    +O(\varphi_k)+O(\|d_k^u\|^{\tau_1})\\
    &=P_k(\nabla^2F(u_k)\!-\!\mathcal{H}_k^u)d_k^u\!+\!O(\|\nabla^2F(u_*)\!-\!\nabla^2F(u_k)\|\cdot \|d_k^u\|)+\\
  &\ \ \ \    O(\varphi_k)+O(\|d_k^u\|^{\tau_1}).
  \end{array}
\end{equation*}
Moreover, since $u_{k+1}=(x_{k+1}, y_{k+1})=(\tilde{x}_{k+1},
\tilde{y}_{k+1})$, the definition (\ref{Hk^}) of $\mathcal{B}_k$,
relations (\ref{B_kx}) and (\ref{setseq}), one can obtain that $\mathcal{B}_k(u_{k+1}-u_*)=0$. This, together with the relationship above, provides
\begin{equation*}
\begin{array}{ll}
  &\left(
  \begin{array}{cc}
      P_k\nabla^2F(u_*)&\mathcal{B}_k^\top \\
      \mathcal{B}_k &0
    \end{array}
    \right)
    \left(
    \begin{array}{c}
      u_{k+1}-u_* \\
      0
    \end{array}
    \right)\\
 & \ \ \ =P_k(\nabla^2F(u_k)-\mathcal{H}_k^u)d_k^u+O(\|\nabla^2F(u_k)-\nabla^2F(u_*)\|\cdot \|d_k^u\|)+\\
  &\ \ \ \ \ \   O(\varphi_k)+O(\|d_k^u\|^{\tau_1}).
\end{array}
\end{equation*}
This, together with Lemma \ref{lemma6}, implies that
\begin{equation}\label{N15}
\begin{array}{ll}
\|u_{k+1}-u_*\|&\!=\!O(\|P_k(\mathcal{H}_k^u\!-\!\nabla^2F(u_k))d_k^u\|)\!+\!O(\|
\nabla^2F(u_k)\!-\!\nabla^2F(u_*)\|\cdot \|d_k^u\|)+\\
& \ \  O(\varphi_k)+O(\|d_k^u\|^{\tau_1}).
\end{array}
\end{equation}
Now, in view of $\nabla^2_{uu}\mathcal{L}_\beta(w_k)=\nabla^2F(u_k)$,
$\tau_1>1$ and $\nabla^2F(u_k)\rightarrow \nabla^2F(u_*)$
(since $F(u)\in {\rm C}^2$), from Assumption \ref{assumption4},
relationships (\ref{F1b})and (\ref{N15}), one has
\begin{equation}\label{N50}
\|u_{k+1}-u_*\|=o(\|u_k-u_*\|)+o(\|d_k^u\|).
\end{equation}
Therefore, in view of $u_{k+1}=u_k+d^u_k$, by Lemma \ref{lemmaN5.4} (ii), the
relationship (\ref{N50}) above implies that $\|u_{k+1}-u_*\|=o(\|u_k-u_*\|)$.
The proof is complete. \hfill$\Box$

The quadratic rate of convergence of the CFD-SQO$_1$-eq-free method is
further analyzed below when the matrics $\mathcal{H}_k^u$ and
$\nabla^2_{uu}\mathcal{L}_\beta(w_k)=\nabla^2 F(u_k)$ have a higher
approximation accuracy.

\begin{assumption}\label{assumption5}
For problem {\rm (\ref{0.1})-eq-free} and the {\rm CFD-SQO$_c$-eq-free}
method, we assume that the following conditions hold.

{\bf(i)} Strengthen Assumption {\rm \ref{assumption3} (i)} to:
the objective functions $f(x)$ and $\theta(y)$ are third-order continuous
differentiable in the neighborhood of the feasible set $\mathcal{F}$;

{\bf(ii)} Same as Assumption {\rm \ref{assumption3} (ii)};

{\bf(iii)} Strengthen Assumption {\rm \ref{assumption4}} to:
the matrix sequence $\{\mathcal{H}^u_k\}$ generated by the method
satisfies one of the following two conditions:

\begin{equation}\label{tiaojian1+}
\|(\mathcal{H}^u_k-\nabla^2 F(u_k))d^u_k\|=O(\|d^u_k\|^2),
\end{equation}
\begin{subequations}\label{tiaojian2+}
  \begin{numcases}{}
      \|P_k(\mathcal{H}^u_k-\nabla^2F(u_k))d^u_k\|=O(\|d^u_k\|^2),\label{tiaojian2a+}\\
      (\alpha^x_*)_{I_*^l}>0, (\gamma^x_*)_{I_*^v}>0, (\alpha^y_*)_{J_*^s}>0, (\gamma^y_*)_{J_*^r}>0, (\mu_*)_{S_*}>0.\label{tiaojian2b+}
  \end{numcases}
  \end{subequations}
\end{assumption}

\begin{theorem}\label{theorem5}
Suppose that Assumptions {\rm \ref{assumption1}} and
{\rm \ref{assumption5}} hold, and parameters $\rho\in (0,0.5),
\tau_1\geqslant 2$, then the CFD-SQO$_1$-eq-free method produces
a sequence $\{u_k\}$ of iterative points satisfying
$\|u_{k+1}-u_*\|=O(\|u_k-u_*\|^2)$, so the {\rm CFD-SQO$_1$-eq-free}
method is quadratically convergent.
\end{theorem}

{\bf Proof}\ \
First, it follows from $F(u)\in {\rm C}^3$ that $\|\nabla^2F(u_k)-\nabla^2F(u_*)\|
=O(\|u_k-u_*\|)$. Second, considering $\tau_1\geqslant 2$ and Assumption \ref{assumption5} (iii),
one has from (\ref{N15}) and (\ref{F1b}) that
\begin{equation}\label{N60}
    \|u_{k+1}-u_*\|
    =O(\|u_k-u_*\|^2)+O(\|u_k-u_*\|\cdot \|d^u_k\|)+O(\|d_k^u\|^2).
\end{equation}
On the other hand, by Theorem \ref{theorem4}, the sequence
$\{u_k\}$ converges superlinearly to $u_*$.
So, by Lemma \ref{lemmaN5.4} (iii), 
one has $\|u_{k+1}-u_k\|\sim \|u_k-u_*\|$,
thus $\|d^u_k\|=\|u_{k+1}-u_k\|\sim \|u_k-u_*\|$. This, together with
relationship (\ref{N60}), shows immediately that  $\|u_{k+1}-u_*\|
=O(\|u_k-u_*\|^2)$. The proof is complete. \hfill$\Box$

\section{Numerical experiments and applications }

In this section, to test the  numerical effects of Algorithm 1, a class of mathematical examples and a kind of economic power dispatch (EPD) are selected to test the numerical performance of Algorithm 1. The numerical
experiment platform is MATLAB R2016a, OPTI2.28 and IPOPT3.12.9, and
the running environment is Intel(R) Core(TM) i5-8500 CPU 3.00GHz RAM
8 GB, Windows 10 (64bite).

\subsection{Experiment with a class of academic examples}
Consider academic model/example (where $q\geq 5$,
sign$(\cdot)$ is the sign function):
\begin{subequations} \label{5.30}
\begin{eqnarray}
\label{5.30a}
\nonumber & \min & F(x):=
\sum_{i=0}^{q-1}\{(2.3x_{(3i+1)}\!+\!0.0001x_{(3i+1)}^2
\!+\!{\rm sign}(q\!-\!5)(-0.0005x_{(3i+1)}^3\!+\!e^{\sin x_{(3i+1)}}))\!+\!\\
\nonumber  &  &   \ \ \ \ \ \ \ \ \ \ (1.7x_{(3i+2)}
\!+\!0.0001x_{(3i+2)}^2\!+\!{\rm sign}(q-5)(0.0008x_{(3i+2)}^3
\!+\!e^{\cos x_{(3i+2)}}))+\\
& & \ \ \ \ \ \ \ \ \ \ (2.2x_{(3i+3)}
\!+\!0.00015x_{(3i+3)}^2\!+\!{\rm sign}(q-5)(0.001x_{(3i+3)}^3
\!+\!e^{\cos x_{(3i+3)}}))\}\label{5.30a}\\
&{\rm s.t.} &   x_1+x_2+x_3 \geq  60,
~x_4 + x_5+ x_6 \geq  50,
~x_7 + x_8+ x_9 \geq 70,\label{5.30b} \label{5.30b}\\
& &  x_{10} + x_{11}+x_{12} \geq  85,
~x_{13} + x_{14}+ x_{15} \geq  100, \label{5.30c} \\
& &  x_{(3i-2)} + x_{(3i-1)}+x_{3i} \geq  100+5(i-4), \ i=6,\dots,q-1, \label{5.30d}\\
& &  -7\leq x_{(3i+1)} - x_{(3i-2)}\leq 6,
~-7\leq x_{(3i+2)} - x_{(3i-1)}\leq 7, \ i=1,\dots,q-1,\label{5.30e}\\
& &  -7\leq x_{(3i+3)} - x_{(3i)}\leq 6,\ i=1,\dots,q-1,
~v_i\leq x_i\leq u_i, \ i=1,\dots,3q, \label{5.30f}
\end{eqnarray}
\end{subequations}
where,
\begin{equation*} 
\label{5.40}
\left\{
\begin{array}{ll}
v_1=8,\ v_2=43,\ v_3=3,\ v_i=0, \ i=4,\dots, 3q; u_1=21,\ u_2=57,\ u_3=16,\\
u_4=u_7=u_{10}=u_{13}=90,\ u_5=u_8=u_{11}=u_{14}=120, u_6=u_9=u_{12}=u_{15}=60,\\
u_{(3i-2)}=90+3i,\ u_{(3i-1)}=120+6i,\ u_{3i}= 60+i, \
 \ i=6,\dots, q.
 \end{array}
\right.
\end{equation*}
The above example was first introduced by Jian et al. \cite{JLY},
which is expanded from the example HS118 in \cite{SXLZ}.
When $q=5$, the model (\ref{5.30}) reduces to HS118, which has an
a known optimal solution and known optimal value, i.e.,
\begin{equation*}
x^*=(8,49,3,1,56,0,1,63,6,3,70,12,5,77,18),\ F(x^*)=664.82045.
\end{equation*}
Therefore, the numerical effect of a tested method can be
examined by solving (\ref{5.30}) with $q=5$.
Divide the variables  into
\begin{equation*}
  \begin{array}{ll}
    \tilde{x}=(x_1,x_4,x_7,\dots,x_{(3q-2)},\ x_2,x_5,x_8,\dots, x_{(3q-1)}),\
    \tilde{y}=(x_3,x_6,x_9,\dots,x_{3q}).
 \end{array}
\end{equation*}
Then model (\ref{5.30}) can be regarded as a two-block structural model
(\ref{0.1}) without the equality constraint (\ref{0.1b}), i.e., $m_1=0$,
constraints (\ref{5.30b})-(\ref{5.30d}) and (\ref{5.30e})-(\ref{5.30f}) correspond to the inequality
constraints (\ref{0.1c}) and (\ref{0.1d})-(\ref{0.1e}), respectively.
The scale of (\ref{5.30}) is $(n_1, n_2; m_1, m_2; l_1, l_2)=(2q, q; 0, q; 4q-2, 2q-1)$.

In the numerical experiments, Algorithm 1 is compared with two closely related algorithms, including the PRS-SQO-DSM \cite{JZY}, the corresponding partially feasible SQO (PF-SQO)
method without the distributed technique (i.e., skipping Steps 1
and 2 in Algorithm 1 and going Step 3 directly). Since the
PRS-SQO-DSM \cite{JZY} is only suitable to handle the optimization
with equality and box constraints, an appropriate equivalence
transformation of the model (\ref{5.30}) is required,
and the transformed model, the approach of variable block,
initial points, parameters, termination criterion and accuracy
are all consistent with the original paper \cite{JZY}.

In the experiments of Algorithm 1 and the PF-SQO method solving
example (\ref{5.30}), the initial point $(\tilde{x},\tilde{y})\in
\mathcal{F_+}$ is generated by solving a linprog function, and
parameters $\tau_1=1.01,M=500, M_1=7000, (\rho,\sigma)=(0.45, 0.9)$.
Due to there exists no equality constraint in model (\ref{5.30}), we would rather choose (\ref{error1}) as the termination criterion with $\epsilon=10^{-8}$. In this way, the accuracy of the algorithm can be improved.
Considering that the Hessian matrix of the objective function of example  (\ref{5.30}) is easy to compute and has a certain ``positive
characterization'', we compute the matrices by
\begin{equation}\label{5.10}
{\mathcal{H}}^x_{k}={H}^x_{k}=\nabla^2 f_{k}+\psi(\nabla^2f_{k})I_{2q},
{\mathcal{H}}^y_{k}={H}^y_{k}=\nabla^2 \theta_{k}+\psi(\nabla^2\theta_{k}) I_{q},
\end{equation}
with $\nabla^2 f_{k}:=\nabla^2 f(x_{k}),
\nabla^2\theta_{k}:=\nabla^2\theta(y_{k})$ and
\begin{equation*}
\label{5.10c}
 \psi(H) =
\left\{
\begin{array}{ll}
0, &  {\rm if}\ \gamma(H)>\eta_0;\\
-\gamma(H)+\eta, &  {\rm if}\ |\gamma(H)|\leq \eta_0;\\
-2\gamma(H), &  {\rm if}\ \gamma(H)<-\eta_0,
\end{array}
\right.
\end{equation*}
where $\eta_0\!=\!0.0001$, $\gamma(D)$ denotes the minimum eigenvalue
of the matrix $D$. As a result, from (\ref{Luua}), one has $\tilde{\mathcal{H}}_k^u=\hat{\mathcal{H}}_k^u=\mathcal{H}_k^u={\rm diag}(H_k^x,\ H_k^y)$. Furthermore, it is clear that the matrices generated above always satisfy Assumption \ref{assumption1}.

The parameter $c$ in Algorithm 1 are chosen as $0, 0.5$ and
$1$, respectively. $C_{t}$ and $F_*$ denote the CPU (second) and the
(approximate) optimal value at the termination
solution, respectively; ${\rm N_{it}}$ and ${\rm N_{Vit}}$ denote
the number of iterations and the number of iterations of
QO splitting (i.e., requirement (\ref{0.19}) is satisfied), respectively,
and  RA denotes the ratio of ${\rm N _{Vit}}$ to ${\rm N_{it}}$;
the feasibility measure of the termination solution
$(\tilde{x}_*,\tilde{y}_*)$ is
\begin{equation*}
\begin{array}{c}
\varphi_{\rm eq}:=\|A\tilde{x}_*+B\tilde{y}_*-b\|_{\infty};\\
\end{array}
\end{equation*}
RE$_{{\rm F}}$ denotes the relative error of $F_*$ generated by the corresponding
method and the PF-SQO method, such as
\[
{\rm RE_F}=\frac{F_* \hbox{ generated by Algorithm 1}-F_* \hbox{ generated by the PF-SQO}}
{F_* \hbox{ generated by the PF-SQO}}\times 100\%.
\]
Calculate the relative error of $C_{\rm t}$ for the PDF-SQO$_{0}$ method and $C_{\rm t}$ for the PF-SQO method:
\[
{\rm RE_C}=\frac{C_t \hbox{ consumed by the PDF-SQO$_{0}$}-
C_t \hbox{ consumed by the PF-SQO}}
{C_t \hbox{ consumed by the PF-SQO}}\times 100\%.
\]

The main numerical results are summarized in Tables \ref{T521}-\ref{T522}.
Based on the data reported in Tables \ref{T521}-\ref{T522},
the following preliminary conclusions can be drawn.

\begin{itemize}

\item In terms of time consumption and number of iterations,
the PDF-SQO$_{1}$ costs the least in the three cases corresponding $c=0,\ 0.5,\ 1$. In terms of total CPU, the PDF-SQO$_{1}$ saves $44.7\%$ and $28.7\%$ compared to the PDF-SQO$_0$ and the PDF-SQO$_{0.5}$, respectively. In terms of the total iterations,
the PDF-SQO$_{1}$ is $47.9\%$ and $31.3\%$ less than the PDF-SQO$_{0}$ and the PDF-SQO$_{0.5}$, respectively.

\item Observing the iterations ${\rm N_{Vit}}$ of the QO splitting,
it is found that almost all NA is above $80\%$, and most is above $90\%$.
This shows that the QO splitting technique in Algorithm 1 is quite efficient.

\item Compared with the PF-SQO, Algorithm 1 has a significant advantage in CPU.
Even if the worst-performing the PDF-SQO$_{0}$ can still saves $78.9\%$ of
computation time compared to the PF-SQO. Moreover, the relative error RE$_{{\rm F}}$
of the (approximate) optimal solution shows that the optimal solution generated
by Algorithm 1 is of higher quality, which indicates that the distributed
SQO method is numerically superior to the non-distributed SQO method.

\item Compared with the PRS-SQO-DSM, the PDF-SQO$_{1}$ has a significant advantage
in CPU, saving $78.3\%$ over the PRS-SQO-DSM.
For the PRS-SQO-DSM, $\varphi_{eq}$ increases as $q$ increases, which indicates that
the feasibility of the generated (approximate) optimal solution is decreasing.
Since (\ref{5.30}) does not contain the equality constraints,
the feasibility of the (approximate) optimal solution generated by Algorithm 1 for (\ref{5.30}) can always be guaranteed,
which is the reason why $F_\ast$ generated by Algorithm 1 is slightly larger than that generated by
the PRS-SQO-DSM. Therefore, Algorithm 1 is an ideal choice when
the feasibility of the (approximate) optimal solution is required
to be more stringent.

\item Algorithm 1 is quite efficient and accurate for model (\ref{5.30}) with $q=5$.
\end{itemize}

In conclusion, Algorithm 1 can efficiently solve the 12
instances of model (\ref{5.30}), and the resulting solutions have
well feasibility and optimality, and require little CPU. The PDF-SQO$_{1}$ has more obvious advantages than PDF-SQO$_{0}$ and PDF-SQO$_{0.5}$,
which indicates that the introduction of the parameter $c$ is beneficial to the
numerical effectiveness of the method.

\begin{sidewaystable}
\setlength{\abovecaptionskip}{0.2cm}  
\setlength{\belowcaptionskip}{0.5cm}
\caption{The experimental numerical results of Algorithm 1,
PF-SQO and PRS-SQO-DSM in solving model (\ref{5.30}), respectively.}
\label{T521}       
\setlength{\tabcolsep}{0.3mm}
\linespread{1.5}\selectfont
\begin{tabular}{c|ccccc|cccc|cccc}
\hline\noalign{\smallskip}
\multirow{2}*{$q$} & \multicolumn{5}{c|}{ PDF-SQO$_{0}$}& \multicolumn{4}{c|}{ PDF-SQO$_{0.5}$}& \multicolumn{4}{c}{ PDF-SQO$_{1}$} \\
\cmidrule{2-14}
 & C$_{\rm t}$ & ${\rm N_{Vit}/N_{it}=RA}$ & ${\rm F_\ast}$ &  ${\rm RE_F}(\%)$&  ${\rm RE_C}(\%)$ &  C$_{\rm t}$ & ${\rm N_{Vit}/N_{it}=RA}$ & $F_\ast$ & ${\rm RE_F}(\%)$ &C$_{\rm t}$ & ${\rm N_{Vit}/N_{it}^{1}=RA}$ & ${\rm F_\ast}$ &  ${\rm RE_F}(\%)$ \\
\hline\noalign{\smallskip}
5 &0.14&17/19=89\%  &  664.8205& -0.00&  -94.89& 0.09& 11/13=85\%  & 664.8205& -0.00&  0.03& 1/3=33\%  & 664.8205& -0.00 \\
50 &0.86&112/123=91\%  &  -100653.91& -0.08&  -90.72& 0.64& 83/94=88\%  & -100661.75& -0.07&  0.51& 56/69=81\%  & -100681.85& -0.05 \\
100&1.33&121/135=90\%  &  -849990.58& -0.22&  -89.72& 1.03& 92/104=88\%  & -849518.79& -0.27&  0.80& 61/74=82\%  & -849864.13& -0.23 \\
200&2.38&145/159=91\%  &  -8584937.26& -0.19&  -91.58& 1.85& 107/121=88\%  & -8584971.78& -0.19&  1.32& 72/83=87\%  & -8573334.45& -0.32\\
300&4.06&166/184=90\%  &  -36275417.44& -0.06&  -90.98& 2.98& 121/136=89\%  & -36257544.63& -0.11&  2.19& 81/96=84\%  & -36259237.63& -0.10\\
400&5.37&187/210=89\%  &  -104013285.34& -0.01&  -90.35& 4.00& 132/150=88\%  & -103987535.86& -0.04&  3.26& 93/115=81\%  & -104011123.57& -0.02\\
500&5.80&190/197=96\%  &  -237336638.13& -0.66&  -91.41& 4.49& 141/148=95\%  & -237392156.15& -0.64&  3.24& 95/102=93\%  & -237454522.21& -0.61\\
600&7.47&205/216=95\%  &  -474411940.13& -0.15&  -90.59& 6.05& 155/167=93\%  & -474619916.50& -0.11&  3.76& 98/103=95\%  & -469877268.57& -1.10  \\
700&8.50&216/226=96\%  &  -852298479.77& -0.18&  -91.42& 6.43& 159/167=95\%  & -850842509.56& -0.35&  4.61& 108/116=93\%  & -850938370.34& -0.34 \\
800&9.46&213/222=96\%  &  -1420203875.26& -0.21&  -91.24& 8.31& 170/179=95\%  & -1420308577.96& -0.20&  5.44& 111/120=93 \% & -1420417660.43& -0.20 \\
900&11.83&241/250=96\%  &  -2234220685.04& -0.19&  -90.66& 9.03& 180/189=95\%  & -2234311731.30& -0.19&  6.61& 121/132=92\%  & -2236570863.26& -0.09\\
1000&13.18&249/258=97 \%&  -3356197867.71& -0.18&  -90.62& 10.66& 190/201=95 \%& -3359371301.83& -0.08&  7.15& 125/133=94 \%& -3353871312.07& -0.25  \\
\noalign{\smallskip}\hline
Total  & 70.38 &~ 2062/2199=94\% & & & & 54.56&~ 1541/1669=92\%&  & & 38.92& 1022/1146=89\%& \\
\noalign{\smallskip}\hline
\end{tabular}
\end{sidewaystable}

\begin{table}
\setlength{\abovecaptionskip}{0.2cm}  
\setlength{\belowcaptionskip}{0.5cm}
\caption{The experimental numerical results of Algorithm 1,
PF-SQO and PRS-SQO-DSM in solving model (\ref{5.30}), respectively
(continued to Table 1).}
\label{T522}       
\setlength{\tabcolsep}{3.5mm}
\linespread{1.5}\selectfont
\begin{tabular}{c|cc|ccc}
\hline\noalign{\smallskip}
\multirow{2}*{$q$} &\multicolumn{2}{c|}{PF-SQO}& \multicolumn{3}{c}{PRS-SQO-DSM~ \cite{JZY}}\\
\cmidrule{2-6}
 & C$_{\rm t}$ &  ${\rm F_\ast}$ &  ${\rm C_t}$ & $\varphi_{\rm eq}$  & ${\rm F_\ast}$ \\
\hline\noalign{\smallskip}
5 & 2.74 &  664.8205&  33.53& 0.05& 664.4469 \\
50 & 9.26 &  -100731.9991& 5.21& 0.10& -100805.95 \\
100 & 12.97 &  -851848.4624& 7.16& 0.41& -852125.53 \\
200 &28.26 &  -8601024.8229&  7.12& 1.66& -8607951.18 \\
300 & 45.03 &  -36296013.8192& 8.64& 3.79& -36349214.05 \\
400 & 55.64 &  -104028525.8658& 10.74& 6.78& -104253825.26 \\
500 & 67.50 &  -238919860.1155& 8.59& 10.54& -239606919.50 \\
600 & 79.43 &  -475121472.8440&  15.35& 15.20& -476837139.87 \\
700 & 98.96 &  -853814075.5936&  12.30& 20.62& -857511463.20 \\
800 & 108.05 &  -1423211419.9339&  12.87& 26.85& -1430403970.65 \\
900 & 126.73 &  -2238551531.1670&  26.94& 33.99& -2251518625.06 \\
1000& 140.47 &  -3362108667.9320&  31.00& 41.84& -3384003612.77 \\
\noalign{\smallskip}\hline
Total &775.04& &179.45&161.83 \\
\noalign{\smallskip}\hline
\end{tabular}
\end{table}

\subsection{Applications in economic power dispatch}

The economic power dispatch (EPD) model is a  dispatch (generation) scheme
designed to minimize the total generation cost of the power supply
system under the physical and system constraints of the unit, as well as
in the state of determining the start and stop state of the unit;
more details see \cite{31,32a,32b}. A broader class of EPD model can be expressed as
\begin{subequations} \label{5.5}
\begin{eqnarray}
\label{5.5a} &\min &F(P)\!:=\!\sum_{i=1}^{N}\sum_{t=1}^{T}\{a_iP_{i,t}^{3}\!+\!b_{i} P_{i,t}^{2}\!+\!c_i P_{i,t}\!+\!d_i\!+\!e_i\left|\sin(f_i(P_{i,t}\!-\!P_{i,\min}))\right|^\delta\}\}  \\
&{\rm s.t.} & \sum_{i=1}^{N}P_{i,t}=P_{D,t},\ t\in \{1,  \dots, T\},\label{5.5b} \\
&& -D_{i}\leq P_{i,t}-P_{i,t-1}\leq U_{i},\ i\in \{1, \dots, N\},\ t\in \{1, \dots, T\},
 \label{5.5c} \\
&& P_{i,\min}\leq P_{i,t}\leq P_{i,\max},\ i\in \{1, \dots, N\},\ t\in \{1, \dots, T\}.
\label{5.5d}
\end{eqnarray}
\end{subequations}

In model (\ref{5.5}), 
$F(P)$ is the total cost of system operation, $P_{i,t}$ is
the generation capacity of unit $i$ in time period $t$, and $P_{i,0}=0.5P_{i,\max}$,
$N$ is the number of
units in the system, $T$ is the number of time periods for optimal scheduling,
$a_{i}, b_{i}, c_{i}, d_{i}, e_{i}$ and
$f_{i}$ are the associated cost coefficients of unit $i$, and
$e_i\left|\sin(f_i(P_{i,t} - P_{i,\min}))\right|^\delta$ is the cost of the valve-point effect.
In particular, the model ignores the valve point effect when $e_i=0$.
Constraint (\ref{5.5b}) is the power balance constraint, and $P_{D, t}$ is the
grid-wide load for the time period $t$. Constraint (\ref{5.5d}) is the upper and lower
generator output constraint, and $P_{i, \min}, P_{i, \max}$ are the
minimum and maximum output of generator $i$, respectively.
Constraint (\ref{5.5c}) is the climbing constraint of the generator, and $D_i, U_i$ are
the upward and downward climbing rate limits of the generator $i$,
respectively.

The parameter $\delta$ in the cost function (\ref{5.5a})
is usually taken as $1$ or $2$. When $\delta=1$ and $e_i\not=0$,
(\ref{5.5}) is a non-smooth EPD model with a valve-point effect.
The model (\ref{5.5}) with $\delta=2$ was first considered by
Jian et. al in \cite{JOTA}, which can be regarded as an approximate
smoothing of the one with $\delta=1$ or a new exploration of EPD model with valve-point effect. Ref. \cite{JOTA} has shown that model (\ref{5.5}) with $\delta=2$ produces very similar effects to  $\delta=1$.

Taking $N_{1}=[\frac{N}{2}], N_{2}=N-N_{1}$, then the variables
$P_{_{i,t}}$ of (\ref{5.5}) are divided into the group of $x$ and
$y$ as follows.
\begin{equation*}
  \begin{array}{ll}
x=(P_{11},P_{12},\ldots,P_{1T},P_{21},P_{22},\ldots,P_{2T},
\ldots,P_{N_11},P_{N_12},\ldots,P_{N_1T})\in \mathbb{R}^{N_{1}T},\\
y=(P_{(N_1+1)1},P_{(N_1+1)2},\ldots,P_{(N_1+1)T},\ldots,P_{N1},P_{N2},
\ldots,P_{NT})\in \mathbb{R}^{N_{2}T}.
  \end{array}
\end{equation*}
And the model (\ref{5.5}) can be reorganized into a two-block model
(\ref{0.1}), with scale $(n_1, n_2; m_1, m_2; l_1, l_2)=(N_1, N_2;
T, 0; 2N_1T, 2N_2T)$, which can be solved by  Algorithm 1.

In the numerical experiments, we generate 20 test instances by copying 5 units of EPD
problem \cite{31} and $T=24$. The structure is shown in Table \ref{T5-1}. The initial iteration point and multiplier are chosen as
\[
u^0=(x^{0},y^{0})
=(P^0_{i,t}=P_{i,\min},\ i\in \{1, \dots, N\},\ t\in \{1, \dots, T\}),\ \lambda^0={\rm ones}(T,1).
\]
Again, considering the equality constraints in model (\ref{5.5}), we prefer to choose (\ref{error2}) as the termination criterion with $\epsilon=0.005$. In this way, the running time and computing cost of the algorithm will be economical.
The matrices $H_{k}^{x}$ and $H_{k}^{y}$ are generated by $\nabla^{2}f(x_{k})$ and $\nabla^{2}\theta(y_{k})$, respectively, i.e., $H_{k}^{y}=\nabla^{2}f(x_{k})$,\   $H_{k}^{y}=\nabla^{2}\theta(y_{k})$.
It is easy to know that the convergence assumptions of Algorithm 1 are valid for this matrix selection.

Algorithm 1, the MS-SQO \cite{JOTA} method, the
PRS-SQO-DSM \cite{JZY} and OPTI solver \cite{OPTI} are used to solve the 20
ED instances in Table \ref{T5-1}. Notice that the model (\ref{5.5})
does not have the inequality constraint corresponding to (\ref{0.1c}),
so Algorithm 1 is independent of parameter $c$ at this point.

First, we consider the model (\ref{5.5}) ignoring the valve-point
effect (i.e., $e_i=0$), and the parameters in Algorithm 1 are chosen as
$\rho=0.49, \ \xi=0.001,\ \beta=200,\ \sigma=0.8$.
The running results are reported in Table \ref{T5-2}.

In the case where the valve-point effect is ignored,
the 20 ED instances tested in Table \ref{T5-2} are convex optimization
problems due to the data of each unit in the objective function is positive.
Next, we further test the effectiveness of Algorithm 1 on six non-convex EPD
instances with the valve-point effect.
In this experiment, the parameters  are selected as
$\rho=0.45, \xi=0.001, \beta=75, \sigma=0.85$,
and the parameters in the valve-point effect are selected as
$\delta=2,\ e_i=10^5a_i,\ f_i=2b_i/(10^5a_i)$.
The computational
results are reported in Table \ref{T5-4}.

In Tables \ref{T5-2} and
\ref{T5-4}, ${\rm RE}$ denotes the relative error of $F_*$ generated
by the corresponding method and OPTI, such as
\[
{\rm RE}=\frac{F_* \hbox{ generated by Algorithm 1}-F_* \hbox{ generated by OPTI}}
{F_* \hbox{ generated by OPTI}}\times 100\%.
\]
It should be noted that in Tables \ref{T5-2} and \ref{T5-4}, the experimental results of
OPTI, MS-SQO method and PRS-SQO-DSM are obtained directly from the
experimental reports in \cite{JZY}.

Based on the analysis of the experimental results reported in Tables \ref{T5-2} and \ref{T5-4}, we can draw the following conclusions.
\begin{itemize}
\item In terms of time consumption, the OPTI solver generally takes longer
time, more than 1000 seconds for instances Nos. 9 to 20, which is
unaccepted for solving realistic economic scheduling problems; Algorithm 1 outperforms the other three methods in most cases due
to its short running time and high solution quality.

\item In terms of the
generated (approximate) optimal solutions and values, all of the relative error RE is negative, i.e. the optimal value generated by Algorithm 1 is better than the one generated by OPTI solver. Further, the maximum error
of Algorithm 1 is less than $-0.0189\%$, so the solution generated by Algorithm 1 has good feasibility and optimality.
\end{itemize}

\section{Conclusions and Perspectives}
By combining the SQO method and the distributed augment Lagrangian method as well as the feasible direction method, this paper proposes a novel partially feasible distributed sequential quadratic optimization method for a class of two-block smooth large-scale optimization problems with both linear equality and linear inequality constraints. The important theoretical characteristics of the proposed method are systematically analyzed and demonstrated, including global convergence, iterative complexity,
superlinear and quadratic rates of convergence. Finally, based on two types of models, the proposed method is tested and the results are satisfactory.

We believe that along with the idea of this paper, there are still some interesting and meaningful problems worth further studying and exploring.

(i) Further study the linear rate of convergence of the proposed method under the K{\L} condition.

(ii) Extend the proposed method to other kind of multi-block optimization problems, such as multi-block optimization problems associated with model (\ref{0.1}), two-block or multi-block optimization problems with nonlinear equality constraint plus linear inequality constraint.

(iii) Study Gauss-Seidel type PFD-SQO method, in which the multipliers are updated two times at each iteration.


\begin{table}
\setlength{\abovecaptionskip}{0.2cm}  
\setlength{\belowcaptionskip}{0.5cm}
\caption{The structures of 20 mid-scale instances obtained by copying the 5-unit system}
\label{T5-1}       
\begin{tabular}{cccccccccccccc}
\hline\noalign{\smallskip}
\multirow{2}*{No.}& \multicolumn{5}{c}{Unit} & {Total } & \multirow{2}*{No.} & \multicolumn{5}{c}{Unit} & {Total}\\
\cmidrule{2-6}
\cmidrule{9-13}
&1&2&3&4&5&{Units $N$}&&1&2&3&4&5&{Units $N$}\\
\hline\noalign{\smallskip}
1&1 &2 &3 & 2 & 2 &10 & 11 & 20 &24 &27 &20 &19 &110\\
2&3 &3 &3 &3 &3 &15 &12 &22 &26 &29 &22 &21 &120\\
3&4 &4 &4 &4 &4 &20 &13 &26 &30 &30 &22 &22 &130\\
4&5 &6 &7 &7 &5 &30 &14 &30 &33 &32 &25 &30 &150\\
5&5 &10 &10 &5 &10 &40 &15 &34 &37 &36 &29 &34 &170\\
6&8 &11 &12 &9 &10 &50 &16 &36 &39 &38 &30 &37 &180\\
7&10 &14 &16 &15 &15 &70 &17 &40 &44 &41 &34 &41 &200\\
8&13 &18 &18 &13 &18 &80 &18 &44 &48 &45 &38 &45 &220\\
9&12 &20 &25 &20 &13 &90 &19 &48 &52 &48 &40 &52 &240\\
10&18 &22 &25 &18 &17 &100 &20 &50 &54 &50 &42 &54 &250\\
\noalign{\smallskip}\hline
\end{tabular}
\end{table}

\begin{table}[htbp]
\footnotesize 
\setlength{\abovecaptionskip}{0.cm}  
\setlength{\belowcaptionskip}{-0.cm} 
\caption{Numerical test results of four methods in solving 20 ED instances without
valve-point effect}
\label{T5-2}       
\setlength{\tabcolsep}{1mm}
\linespread{1.5}\selectfont
\begin{tabular}{ccc|ccc|ccc|ccc}
\noalign{\smallskip}\hline
\multirow{2}*{No.}& \multicolumn{2}{c|}{OPTI } & \multicolumn{3}{c|}{ Algorithm 1 } & \multicolumn{3}{c|}{ MS-SQO \cite{JOTA} }& \multicolumn{3}{c}{ PRS-SQO-DSM \cite{JZY} }\\
\cmidrule{2-12}
& {C$_{\rm t}$}& {$F_*$}& {C$_{\rm t}$}&{$F_*$} 
&{\rm RE(\%)}& {C$_{\rm t}$}  & {$F_*$} & {${\rm RE(\%)}$}  & {C$_{\rm t}$}& {$F_*$} & {${\rm RE(\%)}$}    \\
\noalign{\smallskip}\hline
1& 1.53 &1243485.2    &       0.72& 1243169.69
&   -0.0254&2.48&1243611.17 & 0.0101 & 1.45 & 1243402.72  & -0.0066  \\
2& 2.68&1833617.3    &       0.96& 1833206.75
&   -0.0224 & 2.66 & 1833808.32 & 0.0104 & 1.55& 1833482.67   & -0.0073 \\
3& 5.83 &2444823.1   &         1.17& 2444306.86  
&   -0。0211&2.94&2445089.70  &0.0109 & 1.62& 2445089.42  & 0.0109\\
4& 28.84&3650316.9   &       1.67& 3649583.38
&   -0.0201 &3.94& 3650940.12  &0.0171 & 1.73& 3650600.21   & 0.0078\\
5& 29.57 &5083735.1  &        2.08& 5081810.86
&   -0.0379 &4.69& 5084391.63 &0.0129 & 2.49& 5084610.9  & 0.0172\\
6& 43.50 &6192152.2  &        2.95& 6190828.63
&   -0.0214 &6.86& 6193185.51 &0.0167 & 3.79& 6193348.75  & 0.0193\\
7& 464.27&8636967.3   &        5.34& 8635338.33
&   -0.0189 &17.11&8638532.86  & 0.0181  & 5.30& 8637012.3   & 0.0005\\
8& 654.49&9973328.8  &        6.42& 9970531.52
&   -0.0280 &20.99&9974943.64  &0.0162  & 4.95& 9973855.19  & 0.0053\\
9& 1004.45&11035233.4 &        7.61& 11032974.66
&  -0.0205 &25.50& 11037322.45 &0.0189& 8.45& 11035931.18   & 0.0063\\
10& 1005.87&12291433.0&        7.91& 12288864.36
&  -0.0209  &30.03& 12293741.41 &0.0188 & 10.41& 12292648.26   & 0.0099\\
11& 1004.15&13513839.1&       9.07& 13511013.22
&  -0.0209&36.58&13516385.17 &0.0188 & 11.39& 13515264.88   & 0.0106 \\
12& 1009.19&14736246.1&        10.22& 14733163.45
&  -0.0209&43.06& 14739033.12  & 0.0189  & 13.01& 14737910.31  & 0.0113\\
13& 1009.17&15975567.9&        11.87& 15972451.32
&  -0.0195&50.93& 15978719.24 &0.0197   & 14.99& 15977941.18 & 0.0149\\
14& 1009.88&18492204.7&        14.60& 18488435.24
&  -0.0204&65.53& 18495641.46 &0.0186  & 18.55& 18494616.78  & 0.0130\\
15& 1015.06&20937025.3&        17.71& 20932730.70
&  -0.0205&90.34& 20940974.18 &0.0189  & 22.48& 20939847.58   & 0.0135\\
16& 1015.32&22197414.8&        19.22& 22192857.81
&  -0.0205&99.25& 22201586.65  &0.0188 & 25.53& 22200467.4   & 0.0138\\
17& 1015.29&24659160.1&        22.71& 24654095.49
&  -0.0205&126.55& 24663850.50 &0.0190  & 30.87& 24662733.12  & 0.0145\\
18& 1014.70&27103981.1&        26.37& 27098395.25
&  -0.0206&165.28& 27109222.63 &0.0193 & 36.56& 27108310.74  & 0.0160\\
19& 1015.93&29641688.2&        32.59& 29635556.80
&  -0.0207&202.02& 29647408.79 &0.0193 & 41.66& 29646844.00   & 0.0174\\
20& 1021.03&30864098.3&        35.82& 30857707.43
&  -0.0207&223.97&30870121.91 &0.0195 & 45.45& 30870148.92   & 0.0196\\
\noalign{\smallskip}\hline
Total & 13370.75&  &237.01 & &   & 1220.7 & &  &   302.68 & &  \\
\noalign{\smallskip}\hline
\end{tabular}
\end{table}

\begin{table}[htbp]
\footnotesize  
\setlength{\abovecaptionskip}{0.cm}  
\setlength{\belowcaptionskip}{-0.cm} 
\caption{Numerical test results of four methods in solving 6 non-convex
ED instances with valve-point effect}
\label{T5-4}       
\setlength{\tabcolsep}{1.3mm}
\linespread{2}\selectfont
\begin{tabular}{ccc|ccc|ccc|ccc}
\noalign{\smallskip}\hline
\multirow{2}*{No.}& \multicolumn{2}{c|}{OPTI } & \multicolumn{3}{c|}{Algorithm 1} & \multicolumn{3}{c|}{ MS-SQO \cite{JOTA}} & \multicolumn{3}{c}{PRS-SQO-DSM \cite{JZY}} \\
\cmidrule{2-12}
& {C$_{\rm t}$}& {$F_*$} & {C$_{\rm t}$}&{$F(P^{\ast})$}  & {${\rm RE(\%)}$}  & {C$_{\rm t}$} &{$F_*$} &  {${\rm RE(\%)}$} & {C$_{\rm t}$}& {$F_*$}  & {${\rm RE(\%)}$}  \\
\noalign{\smallskip}\hline
1 & 3.99& 1243486.63  &    0.69& 1242430.80 &	-0.0849 &  3.15 &1243613.78 & 0.0102 & 1.60 &
1243404.31 & -0.0066 \\
6 & 669.00& 6192158.64  &   3.33& 6187183.24 &	-0.0804 & 6.96& 6193193.34  & 0.0167  &	3.55&
6193353.42  & 0.0193 \\
10 & 1009.04 & 12291445.47  &   9.67& 12281422.07  &	-0.0815 & 22.56 & 12293758.39 & 0.0188  &11.91 &
12292659.53& 0.0099\\
14 & 1024.90 & 18492222.66  &   18.77& 18477663.86& -0.0787 & 45.11& 18497628.33  & 0.0292 &	18.48&
18494618.1  & 0.0130 \\
17 & 1114.85 & 24659184.09  & 28.24& 24638774.51  & -0.0828 & 101.11& 24667160.66  & 0.0323   &	31.08& 24662682.15 & 0.0142 \\
20 & 1133.98& 30864170.56  &  44.16& 30839509.37 &	-0.0799 & 181.82& 30874962.95  & 0.0350  &	51.39&
30870016.43  & 0.0189 \\
\noalign{\smallskip}\hline
Total & 4955.8&  & 104.9& &  -0.4882& 360.7& & 0.1422  & 118.01& & 0.0687\\
\noalign{\smallskip}\hline
\end{tabular}
\end{table}

\end{document}